\documentclass[mathpazo]{jiasc}

\usepackage[utf8]{inputenc}
\usepackage[english]{babel}
\usepackage{epsfig,version,graphicx,fancybox,mathrsfs}
\usepackage{subfigure,color,stmaryrd,multirow}
\usepackage{booktabs,siunitx,enumitem}
\usepackage[ruled,vlined]{algorithm2e}
\usepackage{bm,placeins,url}
\usepackage{hyperref}
\hypersetup{hidelinks}
\DeclareGraphicsExtensions{.eps}
\graphicspath{{./}{SO_RtN_Sep2026/}{SO_RtN_Sep2026_EPS/SO_RtN_Sep2026/}}

\makeatletter
\theoremstyle{plain}
\@addtoreset{equation}{section}   

\renewcommand\thefigure{\@arabic\c@figure}
\makeatother
\newtheorem{thm}{\bf Theorem}

\newtheorem{cor}{\bf Corollary}

\newtheorem{lmm}{\bf Lemma}

\renewenvironment{lemma}{\begin{lmm}}{\end{lmm}}
\theoremstyle{remark}

\newcommand{\Om}{\Omega}

\newcommand{\dd}{\,\mathrm{d}}

\def \epsilon {{\varepsilon}}

\definecolor{bgblue}{rgb}{0.04,0.39,0.54}
\definecolor{lired}{rgb}{0.3, 0.0, 0.0}
\definecolor{ligreen}{rgb}{0.0, 0.3, 0.0}
\definecolor{liblue}{rgb}{0.9, 1.0, 1.0}
\definecolor{gray}{rgb}{0.6, 0.6, 0.6}
\definecolor{sky}{rgb}{0.3, 1.0, 1.0}
\definecolor{bunhong}{rgb}{1.0, 0.3, 1.0}
\definecolor{yellow}{rgb}{0.97, 1, 0.0}
\definecolor{liyellow}{rgb}{0.9, 0.8, 0.0}
\definecolor{cengse}{rgb}{0.00,0.40,0.29}

\renewcommand \wedge \times

\begin{document}

\title[A spectral-element method for the Schr\"odinger operator] {
A spectral-element method for computing many eigenvalues and their
asymptotics of the Schr\"odinger operator with Robin boundary condition }
\author[Weizhu Bao and Fukeng Huang]
{Weizhu Bao\affil{1} and Fukeng Huang\affil{2}\comma\corrauth}
\address{\affilnum{1}\ Department of Mathematics, National University of Singapore,
Singapore 119076, Singapore.\\
\affilnum{2}\ School of Mathematical Sciences, Eastern Institute of Technology,
Ningbo, Zhejiang 315200, P. R. China.}
\emails{{\tt matbaowz@nus.edu.sg} (W. Bao), {\tt fkhuang@eitech.edu.cn} (F. Huang)}

\begin{abstract}
We propose spectral and spectral-element methods for accurately computing many
eigenvalues of the Schr\"odinger operator with Robin boundary condition on simple and complex geometries,
respectively. Due to their spectral-type accuracy in approximating those eigenfunctions corresponding to high-index eigenvalues, which are usually highly oscillatory, the proposed approaches have excellent resolution
in computing thousands of eigenvalues accurately and efficiently with the resolution property that the number of eigenvalues with reasonable accuracy is proportional to the number of degrees of freedom.
Using a standard HPC node with 128 GB of memory,
we can obtain numerically more than 5,000 reliable eigenvalues
with relative errors below $10^{-8}$ for different complex domains in two dimensions (2D).
Based on these computed eigenvalues, we first confirm some theoretical results on
the Robin-to-Neumann (RtN) gaps of the Laplacian operator
in 2D, which were recently studied by Rudnick \textit{et al.} [\textit{Comm. Math. Phys.}
388 (2021)]. Then we systematically study the RtN gaps
of the Schr\"odinger operator and their convergence rates.
Based on our extensive numerical results, we formulate
a unified conjecture on the cumulative averages of the RtN gaps of
the Schr\"odinger operator.
\end{abstract}


\maketitle

\section{Introduction}
Consider the following eigenvalue problem of the Schr\"odinger operator with Robin boundary condition, i.e.
find $\lambda \in \mathbb{R}$ and a real-valued function $u:=u({\bf x})\not\equiv0$ such that
\begin{eqnarray}
&&L_{\rm SO}u:=\big[-\Delta +V({\bf x})\big]u({\bf x})=\lambda u({\bf x}),\qquad {\bf x} \in \Omega,\label{SO}\\
&&\frac{\partial u({\bf x})}{\partial {\bf n}}+\sigma u({\bf x})=0,\quad {\bf x} \in \partial \Omega,\label{SORobin}
\end{eqnarray}
where ${\bf x}\in \mathbb{R}^d$ is the spatial coordinate, $\Omega \subset \mathbb{R}^d$ is a connected bounded domain with piecewise smooth boundary, ${\bf n}$ is the unit outward normal vector, $\sigma \in \mathbb{R}$ is a given constant,
and $V:=V({\bf x}) \in L^{\infty}(\Omega)$ is a given real-valued function which
models external fields, trapping effects, impurities or material inhomogeneities.

The Schr\"odinger operator $L_{\rm SO}:=-\Delta +V({\bf x})$ is one of the most
fundamental operators in quantum mechanics, spectral theory and mathematical physics
\cite{bao2019fundamental,berezin1991schrodinger,griffiths2018introduction,teschl2014mathematical2}.
The eigenvalues of the Schr\"odinger operator
represent the admissible energy levels of a quantum system, while the
corresponding eigenfunctions describe the associated stationary states or wave modes.
When $\sigma=0$, the boundary condition \eqref{SORobin} collapses to the homogeneous Neumann boundary condition;
and when $\sigma \to +\infty$, it collapses to the homogeneous Dirichlet boundary condition
$u({\bf x})|_{\partial \Omega}=0$ \cite{bogli2022eigenvalues}. In fact, different boundary conditions encode
different interactions between eigenfunctions and the boundary: (i)
the homogeneous Dirichlet boundary condition corresponds to perfectly confining or impenetrable walls,
(ii) the homogeneous Neumann boundary condition describes zero-flux or reflecting boundaries,
and (iii) the Robin boundary condition models more general boundary interactions,
such as partially reflecting walls, surface potentials and impedance-type
effects \cite{davies1995spectral}. The Schr\"odinger operator \eqref{SO}
with different boundary conditions leads to different spectra and distinct physical interpretations.
Thus understanding the eigenvalues of the Schr\"odinger operator \eqref{SO} with either homogeneous Dirichlet
or homogeneous Neumann or Robin boundary condition is important in quantum confinement, wave propagation,
resonant modes, spectral geometry and related problems in mathematical physics.

Without loss of generality, we assume that $V({\bf x})$ is non-negative, i.e. $V({\bf x})\ge0$ for ${\bf x}\in \Omega$. Since the self-adjoint Schr\"odinger operator with the Robin boundary condition has compact resolvent, its spectrum is real and discrete \cite{babuska1991finite}. We can rank (or order) all eigenvalues (counting multiplicities) $\{\lambda_n^\sigma \ | \ n =1, 2, ...\}$
(with the superscript $\sigma$ referring to the constant in the Robin boundary condition) of
\eqref{SO}-\eqref{SORobin} as
\begin{equation}\label{lbdsg}
\lambda_1^\sigma
\leq\lambda_2^\sigma
\leq\cdots
\leq\lambda_n^\sigma \leq\lambda_{n+1}^\sigma
\leq\cdots,
\end{equation}
where the number of times that an eigenvalue $\lambda$ of \eqref{SO}-\eqref{SORobin}
appears in the above sequence \eqref{lbdsg} is the same as its algebraic multiplicity.
Here we adopt $\{\lambda_n^0 \ | \ n =1, 2, ...\}$ and $\{\lambda_n^\infty \ | \ n =1, 2, ...\}$
to denote the corresponding eigenvalues of the Schr\"odinger operator \eqref{SO} with
homogeneous Neumann ($\sigma=0$ in \eqref{SORobin}) and homogeneous Dirichlet ($\sigma=+\infty$ in \eqref{SORobin})
boundary conditions, respectively. Define the Robin-to-Neumann (RtN) gaps of the Schr\"odinger operator \eqref{SO}-\eqref{SORobin} as
\begin{equation}\label{RtN}
g_n(\sigma):=\lambda_n^\sigma-\lambda_n^0, \qquad n=1,2,\ldots\,,
\end{equation}
and define the cumulative averages of the RtN gaps as
\begin{equation}\label{RtNca}
G_n(\sigma):=\frac{1}{n}\sum_{m=1}^n g_m(\sigma)=\frac{1}{n}\sum_{m=1}^n \left(\lambda_m^\sigma-\lambda_m^0\right),
\qquad n=1,2,\ldots\,.
\end{equation}

When $V({\bf x})\equiv0$ in \eqref{SO}, the Schr\"odinger operator $L_{\rm SO}$ collapses
to the Laplacian operator. In this case, under suitable regularity assumptions on the domain $\Omega$, for any fixed $n\in \mathbb{N}$, one has \cite{bogli2022eigenvalues,ognibene2025asymptotics}
\begin{equation}\label{Robin_parameter_limits}
\lambda_n^\sigma=\lambda_n^0+
\mathcal{O}(|\sigma|), \quad \text{as }\sigma\to0;\qquad
\lambda_n^\sigma=\lambda_n^\infty+\mathcal{O}\left(|\sigma|^{-1}\right),
\ \text{as }\sigma\to +\infty.
\end{equation}
In addition, when the domain $\Omega$ is a simple geometry, all the eigenvalues
of \eqref{SO}-\eqref{SORobin} can be obtained explicitly. Specifically, in one dimension (1D), e.g.
taking $d=1$, $V(x)\equiv 0$ and $\Omega=(0,\pi)$ in \eqref{SO}-\eqref{SORobin}, all the eigenvalues
are given as \cite{gittins2020courant}
\begin{itemize}
\item for the Neumann boundary condition, i.e. $\sigma=0$,
\begin{equation}\label{LONbc}
\lambda_{n}^0=(n-1)^2,\qquad n=1,2,...\,;
\end{equation}

\item for the Dirichlet boundary condition, i.e. $\sigma=+\infty$,
\begin{equation}\label{LODbc}
\lambda_{n}^{\infty}=n^2,\qquad n=1,2,...\,;
\end{equation}

\item for the Robin boundary condition, i.e. $\sigma\ne 0$,
\begin{equation}
\lambda_n^\sigma
=
\begin{cases}
\displaystyle
\frac{2}{\pi}\sigma
+\mathcal{O}(\sigma^2),
& n=1,\quad \sigma\to0^+,\\[2mm]
\displaystyle
(n-1)^2+\frac{4}{\pi}\sigma
+\mathcal{O}(\sigma^2),
& n\geq2,\quad \sigma\to0^+,\\[2mm]
\displaystyle
n^2-\frac{4n^2}{\pi\sigma}
+\mathcal{O}(\sigma^{-2}),
& n\geq1,\quad \sigma\to+\infty.
\end{cases}
\label{LORbc}
\end{equation}
For \(\sigma>0\), the eigenvalues can be written as \(\lambda_n^\sigma=\kappa_n^2(\sigma)/\pi^2\), where \(\kappa_n(\sigma)\) is the unique solution in \(((n-1)\pi,n\pi)\) of
\[
\kappa_n\tan(\kappa_n/2)=\sigma\pi
\]
when \(n\) is odd, and of
\[
\kappa_n+\sigma\pi\tan(\kappa_n/2)=0
\]
when \(n\) is even. For \(\sigma<0\), the first eigenvalue is negative. The second eigenvalue is negative if and only if \(\sigma<-2/\pi\), and is zero at \(\sigma=-2/\pi\). The small-\(\sigma\) expansions in  \eqref{LORbc} hold as \(\sigma\to0\) from either side.
\end{itemize}
Combining \eqref{LONbc}, \eqref{LODbc} and \eqref{LORbc}, \eqref{Robin_parameter_limits} is confirmed immediately
in this special case. In addition, for any $\sigma\in \mathbb{R}$, we can get
\begin{equation}\label{RtN1D1}
\begin{split}
&\lim_{n\to\infty} g_n(\sigma)=\lim_{n\to\infty} (\lambda_n^\sigma-\lambda_n^0)=\frac{4}{\pi}\sigma, \\
&\lim_{n\to\infty} G_n(\sigma)=\lim_{n\to\infty}\frac{1}{n}\sum_{m=1}^n \left(\lambda_m^\sigma-\lambda_m^0\right) =\frac{4}{\pi}\sigma.
\end{split}
\end{equation}
In two dimensions (2D), e.g. taking $d=2$, $V({\bf x})\equiv 0$ and $\Omega\subset\mathbb{R}^2$
to be a connected bounded domain with piecewise smooth boundary in \eqref{SO}-\eqref{SORobin},
Rudnick et al. \cite{rudnick2021differences} studied
the RtN gaps and their cumulative averages. Based on analytical and numerical results for rectangles, the
hemisphere, and the equilateral triangle \cite{rudnick2021rectangles,rudnick2022triangle,rudnick2022hemisphere},
for each fixed $\sigma>0$, they proved the following result
\cite{rudnick2021differences}
\begin{equation}\label{RtNgap}
\lim_{n\to\infty} G_n(\sigma)=\lim_{n\to\infty}\frac{1}{n}\sum_{m=1}^n \left(\lambda_m^\sigma-\lambda_m^0\right)
=\frac{2 |\partial \Omega|}{|\Omega|}\sigma,
\end{equation}
where $|\Omega|$ and $|\partial\Omega|$ denote the area of $\Omega$ and the length of $\partial\Omega$,
respectively. The above formula \eqref{RtNgap} reveals an interesting connection
between an algebraic spectral quantity, i.e. the cumulative averages of the RtN gaps,
and the geometric quantities of the domain $\Omega$, i.e., its area and
boundary length.

In general domains, especially for the Schr\"odinger operator with nontrivial
potentials, exact eigenvalues are rarely available, thus efficient and accurate as well as high resolution
numerical approximation is necessary and important as well as challenging. In fact,
the finite element method provides a mature and
geometrically flexible framework, supported by general Galerkin spectral
approximation theory \cite{babuvska1989,babuska1991finite,boffi2010finite}.
Most concrete algorithms and analyses have been developed for computing the eigenvalues of the
Schr\"odinger operator with homogeneous Dirichlet boundary condition, including adaptive eigensolvers
\cite{carstensen2012adaptive,chen2011finite,dai2008convergence,
yang2021eigenfunction}, two-grid schemes
\cite{hu2011acceleration,xu2001,yang2011two}, multilevel correction methods
\cite{lin2015multi}, and rigorous error estimates and eigenvalue bounds
\cite{carstensen2014guaranteed,carstensen2020skeletal,larson2000posteriori,zhai2019lower}.
For computing the eigenvalues of the
Schr\"odinger operator with homogeneous Neumann condition,
finite element error analysis on curved domains was given in
\cite{hernandez2003neumann}. For computing the eigenvalues of the
Schr\"odinger operator with Robin boundary condition, far fewer results are available in the literature
\cite{lopezyela2017finite,vanmaele1995operator}.
Recently, machine learning methods have also been developed for computing eigenvalues of Schr\"odinger and fractional Schr\"odinger operators \cite{guo2024deep,guo2026generalization}.
Specifically, all these existing
works focus mainly on algorithmic and theoretical analysis, or on a relatively small number of low-index
eigenpairs.

A classical question about the full spectrum is whether the leading term in Weyl's law
bounds each eigenvalue of the Dirichlet Laplacian from below, as predicted by
P\'olya's conjecture \cite{polya1961vibrating}. This bound holds for Euclidean balls
\cite{filonov2023balls} and planar annuli \cite{filonov2026annuli}.
On general bounded domains, the sharp lower bound for the sum of the first $n$
eigenvalues \cite{li1983schrodinger} yields a lower bound for the $n$th eigenvalue.
The constant in this individual bound has been improved using energy orthogonality
\cite{wang2026toward} and sharpened further by coupling Fourier projections
\cite{wang2026improved}. On bounded Lipschitz domains, estimates for the remainder
in Weyl's law yield the conjectured bound up to an $\varepsilon$-loss above an
explicit spectral threshold \cite{jiang2026polya}.

For the Robin Schr\"odinger operator \eqref{SO}--\eqref{SORobin}, we are interested
both in the growth of $\lambda_n^\sigma$ described by Weyl's law and in the differences
$\lambda_n^\sigma-\lambda_n^0$ that define the RtN gaps \eqref{RtN}. To investigate
these gaps and their cumulative averages \eqref{RtNca} numerically, we need accurate
approximations to thousands of Robin and Neumann eigenvalues.
To achieve this very challenging computational task,
under a fixed number of degrees of freedom (DOFs), low-order finite element or finite difference
discretizations typically resolve only a limited portion of the spectrum
\cite{bao2020jacobi,sauter2010hp,zhang2015many}, and thus they cannot handle it
successfully. When the domain $\Omega$ is a simple geometry, e.g. a rectangle or a disc in 2D, and
a cube or a ball in 3D, spectral methods can provide a much larger
range of accurate high-index eigenvalues with the resolution of the number
of computed eigenvalues being proportional to the degrees of freedom of the corresponding
discretization \cite{bao2020jacobi,hashemi2022least,li2017efficient,Shao2026,lin2025highly,
weideman1988eigenvalues,zhang2015many}. In contrast, when the domain $\Omega$ is complicated,
spectral methods cannot be extended directly. As we know, spectral-element methods
can address this geometric limitation by combining
high-order spectral approximation with element-wise decomposition and
geometric mappings \cite{deville2002high,han2013class,patera1984spectral,shen2011spectral,yang2016seamless}.
Representative applications to eigenvalue problems include the Schr\"odinger
operator with singular potentials and homogeneous Dirichlet boundary condition
\cite{li2017efficient,wang2022spectral}, and the Stokes
eigenvalue problems with homogeneous Dirichlet boundary condition
\cite{shan2017triangular}.

To the best of our knowledge, efficient and accurate
numerical methods for computing many eigenvalues of the Schr\"odinger operator
with Robin boundary condition are very limited in the literature.
The purpose of this paper is: (i) to develop spectral and spectral-element methods for accurately computing many
eigenvalues of the Schr\"odinger operator with Robin boundary condition on simple and complex geometries,
respectively; and (ii) to adopt the proposed numerical method for obtaining
numerically more than 5,000 reliable eigenvalues
with relative errors below $10^{-8}$ for different complex domains in two dimensions (2D)
using a standard HPC node with 128 GB of memory and thus to study numerically
the Robin-to-Neumann gaps and their asymptotics. Specifically, based on our extensive numerical results and observations, we speculate the following:

\medskip
{\bf RtN Gap Conjecture} (Robin-to-Neumann gaps of the Schr\"odinger operator \eqref{SO}-\eqref{SORobin})
Assume that $\Omega\subset\mathbb{R}^d$ is a connected bounded domain with
piecewise smooth boundary and $V\in L^\infty(\Omega)$.
For any $\sigma\in\mathbb{R}$, we have
\begin{equation}\label{RtNgapdim}
\lim_{n\to\infty} G_n(\sigma)
=\frac{2 |\partial \Omega|}{|\Omega|}\sigma, \qquad G_n(\sigma)=\frac{2 |\partial \Omega|}{|\Omega|}\sigma
+O(n^{-1/d}), \quad n\gg1.
\end{equation}
In addition, when $d=1$ and $\Omega=(a,b)$, we also have
\begin{equation}\label{RtNgap1d}
\lim_{n\to\infty} g_n(\sigma)=\lim_{n\to\infty}(\lambda_n^\sigma-\lambda_n^0)
=\frac{4}{b-a}\sigma, \qquad \lambda_n^\sigma=\lambda_n^0+\frac{4}{b-a}\sigma
+O(n^{-2}), \quad n\gg1,
\end{equation}
and
\begin{equation}
G_n(\sigma)=\frac{4\sigma}{b-a}-(\sigma^2+\frac{2\sigma}{b-a})n^{-1}+O(n^{-2}),\quad n \gg 1.
\end{equation}

\medskip

The rest of the paper is organized as follows. In Section~2, we introduce the scaling property and the variational formulation of the
Schr\"odinger eigenvalue problem \eqref{SO}-\eqref{SORobin}. In Section~3, we describe the spectral and spectral-element methods for solving the problem \eqref{SO}-\eqref{SORobin} on simple and complex geometries, respectively. Section~4 is devoted to accuracy tests and comparisons of different combinations of the number of sub-elements
and the polynomial degree of basis functions under fixed numbers of DOFs. In Section~5, we investigate the Robin-to-Neumann gaps in 1D, 2D and three dimensions (3D). Finally, some concluding
remarks are given in Section~6.

\section{Scaling property and the variational formulation}

In this section, we present a scaling property of the Schr\"odinger eigenvalue problem \eqref{SO}-\eqref{SORobin}
and its variational (weak) formulation.

\subsection{Scaling property}
Denote
\begin{equation}\label{mapping}
\begin{split}
&R:={\rm diam}(\Omega)=\max_{{\bf x},{\bf z}\in\bar{\Omega}}|{\bf x}-{\bf z}|,\quad {\bf y}:=\frac{{\bf x}}{R},\quad
\tilde \Omega:=\{{\bf y}:={\bf x}/R \ |\ {\bf x}\in \Omega\}, \\
&\tilde V({\bf y}):=R^2V({\bf x})=R^2 V(R{\bf y}), \quad v({\bf y}):=u({\bf x})=u(R{\bf y}), \quad
\tilde \lambda =R^2\lambda, \quad \tilde \sigma=R\sigma.
\end{split}
\end{equation}
Plugging \eqref{mapping} into the
Schr\"odinger eigenvalue problem \eqref{SO}-\eqref{SORobin}, we obtain the following
Schr\"odinger eigenvalue problem defined on a bounded domain $\tilde \Omega$ with unit diameter, i.e.
find $\tilde\lambda \in \mathbb{R}$ and a real-valued function $v:=v({\bf y})\not\equiv0$ such that
\begin{eqnarray}
&&\left[-\Delta +\tilde V({\bf y})\right]v({\bf y})=\tilde \lambda v({\bf y}),
\qquad  {\bf y}\in\tilde \Omega, \label{SOscal}\\
&&\frac{\partial v({\bf y})}{\partial {\bf n}}+\tilde \sigma v({\bf y})=0,\qquad  {\bf y}\in\partial\tilde \Omega.
\label{SOscal234}
\end{eqnarray}

Then it is straightforward to have the following lemma:
\begin{lemma}
Let $\tilde \lambda$ be an eigenvalue of \eqref{SOscal}-\eqref{SOscal234} and $v:=v({\bf y})$ be the corresponding eigenfunction, than $\lambda =\tilde \lambda/R^2$ is an eigenvalue of \eqref{SO}-\eqref{SORobin} and
$u:=u({\bf x})=v({\bf y})=v({\bf x}/R)$ is the corresponding eigenfunction, and the converse also holds.
\end{lemma}

From this lemma, we can easily obtain that: assume that $\tilde\lambda_1^{\tilde\sigma} \leq\tilde \lambda_2^{\tilde\sigma}\leq\cdots
\leq\tilde\lambda_n^{\tilde\sigma} \leq\cdots$ are all eigenvalues of \eqref{SOscal}-\eqref{SOscal234},
then $\lambda_1^\sigma  \leq\lambda_2^\sigma \leq\cdots\leq\lambda_n^\sigma
\leq\cdots$ are all eigenvalues of \eqref{SO}-\eqref{SORobin} with $\lambda=\tilde \lambda/R^2$ and
$\sigma=\tilde \sigma/R$.

\subsection{A variational (weak) formulation} Multiplying a test function $w\in H^1(\Omega)$ in both sides of
\eqref{SO} and then integrating over $\Omega$ and using integration by parts, we obtain a variational formulation
of the problem \eqref{SO}-\eqref{SORobin} as: find $\lambda \in \mathbb{R}$ and a real-valued function $0\ne u\in H^1(\Omega)$ such that
\begin{equation}\label{weak}
a_\sigma(u,w)=\lambda\, b(u,w), \qquad \forall w \in H^1(\Omega),
\end{equation}
where
\begin{equation}\label{bilinear}
a_\sigma(u,w)=\int_{\Omega}(\nabla u \cdot \nabla w+Vuw)\, d{\bf x}+ \sigma\int_{\partial \Omega} uw ds,\quad
 b(u,w)=\int_{\Omega} uw d{\bf x},\quad \forall w \in H^1(\Omega).
\end{equation}

When $\sigma\ge0$, then $a_\sigma(v,v)\ge a_0(v,v)$; and respectively, when $\sigma<0$,
from the trace inequality \cite[Theorem~1.6.6]{brenner2008mathematical},  $\forall\,\varepsilon>0$, we have
\begin{equation}
 \|u\|^2_{L^2(\partial \Omega)}
 \le C \|u\|_{L^2(\Omega)}\|u\|_{H^1(\Omega)}
 \le C_{\varepsilon}\|u\|^2_{L^2(\Omega)}
 +\varepsilon\|\nabla u\|^2_{L^2(\Omega)}.
 \label{eq:trace-eps}
\end{equation}
Since $\varepsilon>0$ can be chosen sufficiently small in \eqref{eq:trace-eps}, the negative boundary
term in $a_\sigma(v,v)$ can be controlled by  $\|\nabla u\|^2_{L^2(\Omega)}$. Assume $V\in L^\infty(\Omega)$,
then the smallest eigenvalue of \eqref{weak} is bounded below and all eigenvalues of \eqref{weak} are discrete.
Specifically, from the Rayleigh quotient \cite{boffi2010finite}, the smallest eigenvalue $\lambda_1^\sigma$ is
defined as
\begin{equation}\label{RQ}
\lambda_1^\sigma=\min_{0\ne v\in H^1(\Omega)}\;\frac{a_\sigma(v,v)}{b(v,v)}=\min_{v\in H^1(\Omega)\; \&\; b(v,v)=\|v\|^2_{L^2(\Omega)}=1}\; a_\sigma(v,v).
\end{equation}
For any $n\ge1$, assume that $u_1({\bf x}), \ldots, u_n({\bf x})$ being the orthonormal eigenfunctions corresponding
to the eigenvalues $\lambda_1^\sigma\le \ldots\le \lambda_n^\sigma$, and define $U_n={\rm span}\{u_1, \ldots, u_n\}$, then
from the Rayleigh quotient \cite{boffi2010finite}, the next eigenvalue $\lambda_{n+1}^\sigma$ is
defined as
\begin{equation}\label{RQn1}
\lambda_{n+1}^{\sigma}=\min_{0\neq v\in U_n^\perp}\; \frac{a_\sigma(v,v)}{b(v,v)}=\min_{E_{n+1}\subset H^1(\Omega)}\;
\max_{0\neq v\in E_{n+1}}\; \frac{a_\sigma(v,v)}{b(v,v)};
\end{equation}
where the dimension of the subspace $E_{n+1}$ is $n+1$.
Combining \eqref{bilinear}, \eqref{RQ} and \eqref{RQn1}, we obtain immediately: (i) when $\sigma>0$, we have
\begin{equation}
\lambda_n^{0} \le \lambda_n^{\sigma} \le \lambda_n^{\infty},\qquad  n=1,2,\ldots\;;
\end{equation}
and respectively, (ii) when $\sigma<0$, we have
\begin{equation}
\lambda_n^{\sigma}\le \lambda_n^{0} \le \lambda_n^{\infty},\qquad  n=1,2,\ldots\,.
\end{equation}

\section{Spectral and spectral-element methods}\label{numerical}
In this section, we briefly introduce the spectral and spectral-element methods
for solving the eigenvalue problem \eqref{SO}-\eqref{SORobin} on simple and complex/irregular
geometries, respectively.
For a more detailed description of the implementation, we refer the reader to
Chapter~4 of \cite{shen2011spectral} and Section~4.5 of \cite{deville2002high}.

\subsection{A Galerkin approximation}
Let $X_N={\rm span}\{\phi_1,\ldots,\phi_{N}\}\subset H^1(\Om)$ be
a conforming finite-dimensional subspace with the degree of freedoms (DOFs) as $N$.
Then a Galerkin approximation of \eqref{weak} reads: find $\lambda_N \in \mathbb{R}$ and a real-valued function $0\ne u_N\in X_N$ such that
\begin{equation}\label{eigenfinite}
  a_{\sigma}(u_N,w_N)=\lambda_N\, b(u_N,w_N), \qquad \forall w_N\in X_N.
\end{equation}
Assume
\begin{equation}\label{uhterm}
u_N=\sum_{j=1}^{N}u_j\phi_j= \Phi^TU\in X_N, \ {\rm with}\ \Phi=(\phi_1,\ldots,\phi_{N})^T, \
U=(u_1,\ldots, u_{N})^T\in {\mathbb R}^{N},
\end{equation}
Plugging \eqref{uhterm} into \eqref{eigenfinite} and taking $w_N=\phi_i$ for $i=1,\ldots,N$ respectively,
we obtain the following generalized matrix eigenvalue problem
\begin{equation}\label{eq:matrix-evp}
  A\, U=\lambda_N\, M\, U,
\end{equation}
where $ A=(a_{ij})\in {\mathbb R}^{N\times N}$ and  $M=(m_{ij})\in {\mathbb R}^{N\times N}$ are two $N\times N$ matrices given by
\begin{equation}\label{aijbij}
  a_{ij}=a_{\sigma}(\phi_j,\phi_i),
  \qquad
  m_{ij}=b(\phi_j,\phi_i), \qquad i,j=1,\ldots, N.
\end{equation}

Then spectral and spectral-element methods can be
designed by properly choosing the basis functions $\{\phi_1,\ldots,\phi_{N}\}$,
efficiently and accurately evaluating the integrals $a_{\sigma}(\phi_j,\phi_i)$ and
$b(\phi_j,\phi_i)$, and properly assembling the two matrices $A$ and
$ M$. Based on different domains $\Omega$, we will present details in
the next two subsections.

\subsection{Spectral methods for simple geometries}

When the domain $\Omega$ is a simple geometry, e.g. an interval in 1D,  a rectangle or a disk in 2D,
a cube or a ball or a cylinder in 3D, spectral methods can be easily designed
for the Galerkin approximation \eqref{eigenfinite}.

For simplicity of notation,
here we only present the spectral method for $\Omega=(-1,1)$ in 1D. In this case,
we adopt the Legendre polynomials to construct the basis in $X_N$. Define
 \begin{equation}\label{Legen}
  \phi_1(x)=\frac{1-x}{2},\quad
  \phi_2(x)=\frac{1+x}{2},\quad
  \phi_i(x)=L_{i-1}(x)-L_{i-3}(x),\ i\geq3, \quad x\in[-1,1],
\end{equation}
where $L_j$ denotes the Legendre polynomial of degree $j\ge0$ \cite{shen2011spectral},
then we can construct a subspace $X_N={\rm span}\{\phi_1,\ldots,\phi_N\}\subset H^1(\Omega)$
with spectral accuracy.
Plugging \eqref{Legen} into \eqref{eigenfinite}, \eqref{bilinear} and \eqref{aijbij},
 we obtain the corresponding
generalized matrix eigenvalue problem \eqref{eq:matrix-evp} with
\begin{equation}
\begin{split}
  &a_{ij}=\int_{-1}^1\left[\phi_j'(x)\phi_i'(x)+V(x)\phi_j(x)\phi_i(x)\right]\dd x+\sigma\left[\phi_j(-1)\phi_i(-1)+\phi_j(1)\phi_i(1)\right],\\
  &m_{ij}=\int_{-1}^1\phi_j(x)\phi_i(x)\dd x, \qquad i,j=1,\ldots, N.
\end{split}
\end{equation}
In practical computations, the above integrals will be further approximated by numerical
quadrature (pseudospectral method) \cite{shen2011spectral}.

The above spectral method for discretizing the eigenvalue problems \eqref{weak} can be extended straightforwardly to a rectangle or a disk in 2D, and
a cube or a ball or a cylinder in 3D, by adopting the standard spectral basis functions
in the literature \cite{shen2011spectral}.
\subsection{Spectral-element methods for complex geometries}
For complex geometries, we adopt conforming spectral-element methods.
For simplicity, we describe the construction for a planar domain in 2D \cite{yang2016seamless};
the extension to 3D follows the same tensor-product framework
\cite{deville2002high}.

\begin{figure}[!htbp]
 \begin{center}
  \subfigure{ \includegraphics[scale=.4]{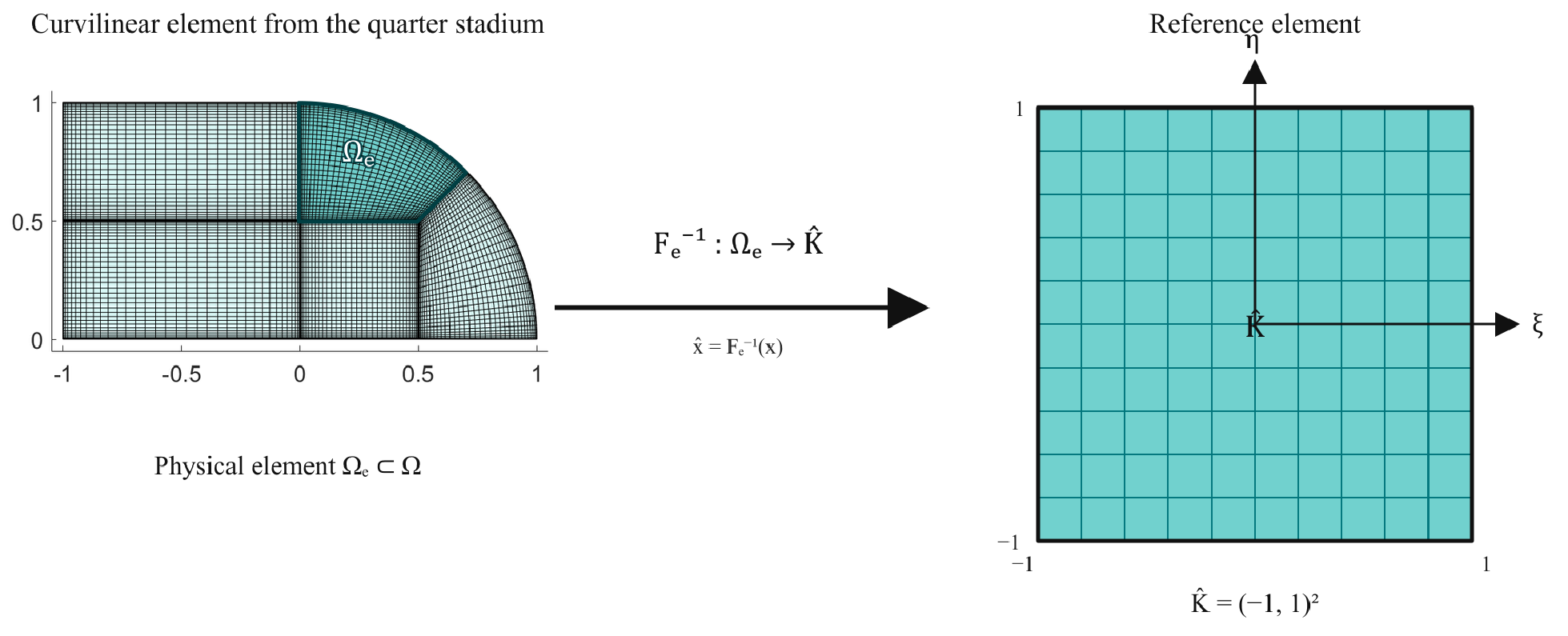}}
  \end{center}
 \caption{Sketch of coordinate transformation from physical domain $\Om_e$ to computational domain $\widehat K$.}\label{fig:SPE}
 \end{figure}

We first decompose the domain into $N_E$ non-overlapping
curvilinear quadrilateral elements as
$\overline\Om=\bigcup_{e=1}^{N_E}\overline\Om_e$,
and introduce smooth element maps (cf. Fig. \ref{fig:SPE})
\[
  F_e:\widehat K=(-1,1)^2\longrightarrow \Om_e.
\]
The four edges of each element are described by smooth
parameterizations $\boldsymbol\gamma_{e,r}(t)$, $t\in[-1,1]$,
$r=1,\ldots,4$, with matching endpoints.
We construct $F_e$ by Gordon--Hall blending of these edge curves,
using compatible parameterizations on shared edges.
The curved-element construction follows
Yang et al.~\cite{yang2016seamless}.

On the reference square $\widehat K$, we can use tensor-product Gauss--Lobatto--Legendre (GLL) nodes and
the associated Lagrange basis with the degree of polynomials in each direction up to $N_p$.
The local trial functions on $\Om_e$ are
obtained by composition with $F_e^{-1}$, while the tensorial GLL grid
is mapped to each physical element by $F_e$.
Let $\widehat u=u\circ F_e$ and $\widehat w=w\circ F_e$.  The local volume form is
\begin{equation}
a_e^{\mathrm{vol}}(u,w)
=\int_{\widehat K}
(\nabla_{\widehat {\bf x}}\widehat u)^T G_e
(\nabla_{\widehat {\bf x}}\widehat w)
\dd\widehat {\bf x}
\quad+
\int_{\widehat K}
\widehat V_e(\widehat {\bf x})
\widehat u(\widehat {\bf x})\widehat w(\widehat {\bf x})
\,|\det(D\,F_e(\widehat {\bf x}))|\dd\widehat {\bf x},
\end{equation}
where
\[G_e=
  |\det (D\,F_e)|\,
  (D\,F_e)^{-1}((D\,F_e)^T)^{-1},
  \qquad
  \widehat V_e=V\circ F_e .
\]
The local mass is
\[
  m_e(u,w)
  =
  \int_{\widehat K}
  \widehat u(\widehat {\bf x})\widehat w(\widehat {\bf x})
  \,|\det(DF_e(\widehat {\bf x}))|\dd\widehat {\bf x}.
\]
If an edge $\gamma$ of $\Om_e$ lies on the physical boundary $\partial \Omega$, the Robin term
adds
\[
  \sigma
  \int_{\widehat\gamma}
  \widehat u\,\widehat w\,
  |\partial_\tau F_e|\dd\widehat s,
\]
where $\widehat\gamma$ is the corresponding reference edge and
$\partial_\tau$ denotes differentiation along that reference edge.
The volume and boundary integrals are evaluated by GLL quadrature,
following \cite{yang2016seamless}.
The local matrices, including the Jacobian and metric factors, are
assembled by identifying the nodal unknowns on shared vertices and
edges to enforce continuity. Robin contributions are added only on
$\partial\Omega$. This gives the generalized matrix eigenvalue problem
\eqref{eq:matrix-evp}, which we solve using MATLAB's \texttt{eigs}.

\section{Accuracy test}
In this section, we report accuracy of the spectral and spectral-element methods
proposed in the previous section for computing the Schr\"{o}dinger
eigenvalue problem \eqref{SO}-\eqref{SORobin}. Let $N$ be the degree of freedoms (DOFs)
of an approximation, and $\lambda_{n,N}^{\sigma} (n=1,2,\ldots, N)$ be all the eigenvalues
of the numerical approximation \eqref{eigenfinite}, satisfying (or ranked as)
\[\lambda_{1,N}^{\sigma}\le \lambda_{2,N}^{\sigma}\le \ldots \le \lambda_{N,N}^{\sigma}.
\]
We adopt the relative errors which are defined as
\begin{equation}
{\rm Relative\ error}\ (\lambda_{n}^{\sigma}):= \frac{|\lambda_{n}^{\sigma}-\lambda_{n,N}^{\sigma}|}{|\lambda_{n}^{\sigma}|}, \qquad n=1,2,\ldots,N,
\end{equation}
if $\lambda_{n}^{\sigma}\neq0$, and otherwise we use the absolute error as $|0-\lambda_{n,N}^{\sigma}|$
if $\lambda_{n}^{\sigma}=0$. In the case that the eigenvalues $\lambda_{n}^{\sigma}$ are not given analytically,
the numerical `exact' eigenvalues are obtained by choosing a sufficient large $N=N_{\rm ex}$.

\subsection{Accuracy in 1D} We take $d=1$ and  $\Omega=(0, \pi)$ in \eqref{SO}-\eqref{SORobin}.
When $V(x)\equiv 0$, all the eigenvalues are given analytically in \eqref{LONbc}, \eqref{LODbc}
and \eqref{LORbc} for homogeneous Neumann boundary condition, homogeneous Dirichlet boundary
condition and Robin boundary condition, respectively \cite{gittins2020courant}. We adopt the
spectral method presented in the previous section to solve this problem numerically for obtaining
the corresponding eigenvalues.

\begin{figure}[!htbp]
 \centering
  \subfigure{ \includegraphics[scale=.3]{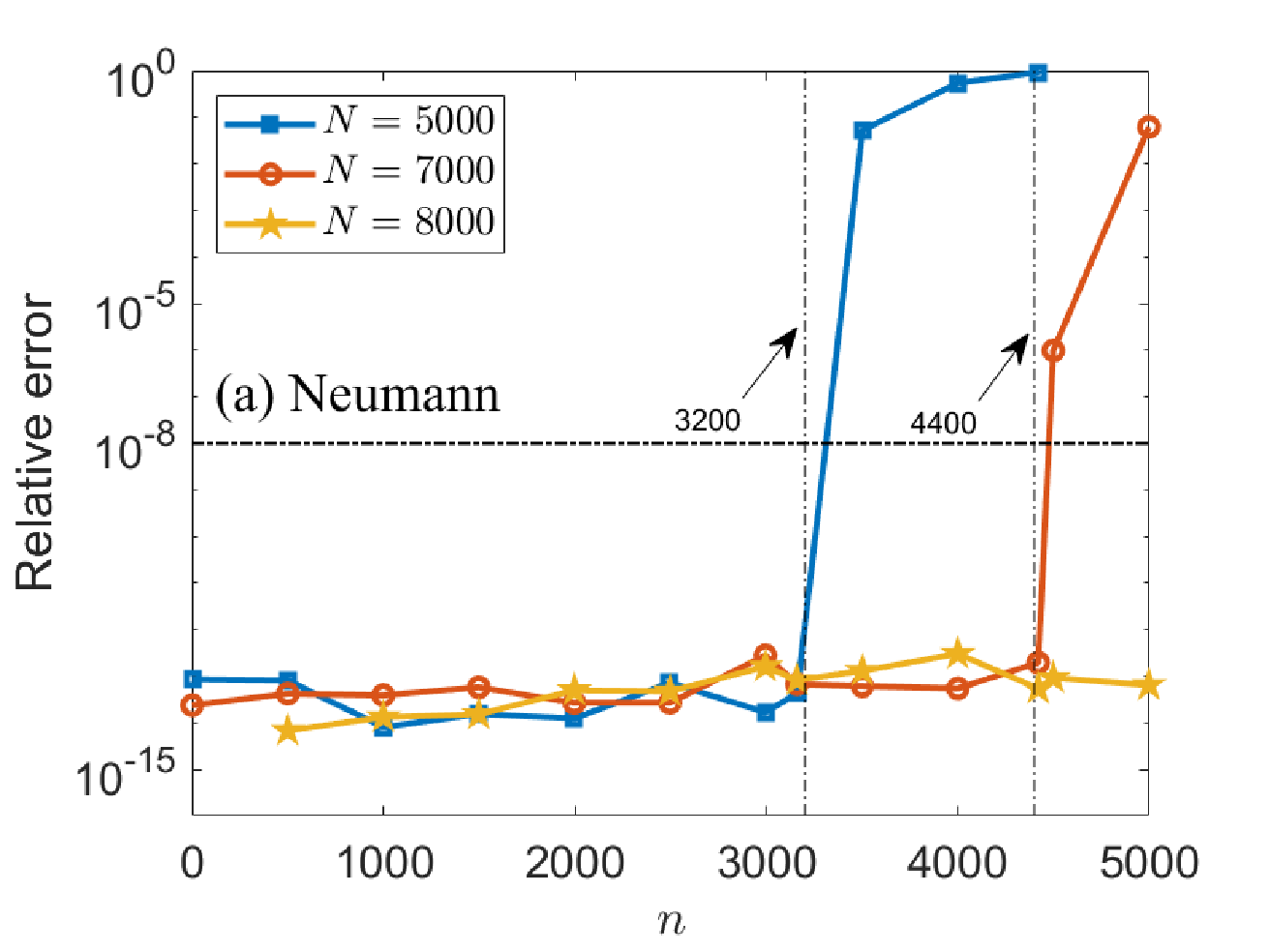}}
  \subfigure{ \includegraphics[scale=.3]{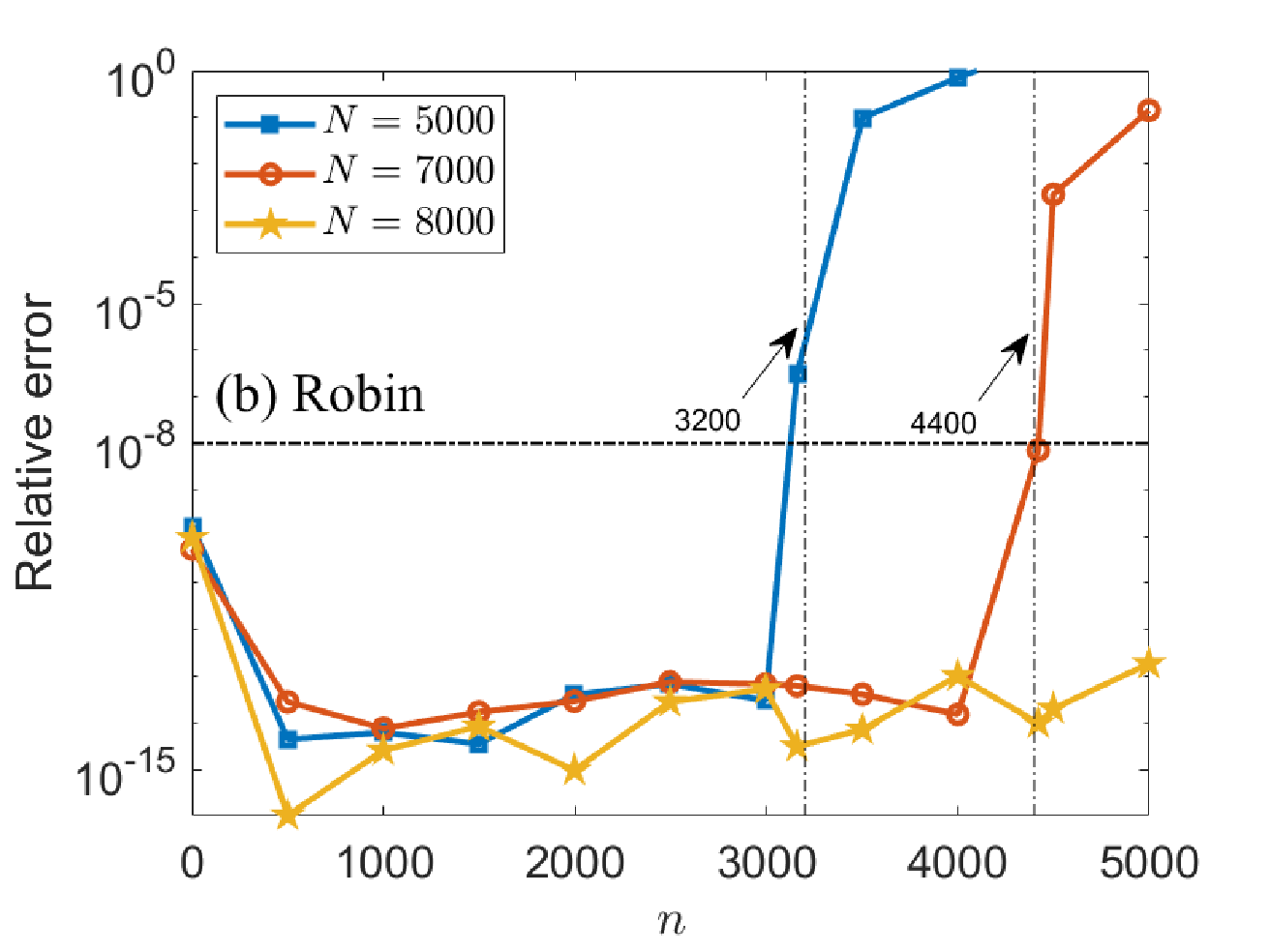}}
  \caption{Accuracy test in 1D with
$\Omega=(0,\pi)$ and  $V(x)\equiv0$ in \eqref{SO}-\eqref{SORobin} for different $N$ and $\sigma$:
(a) Neumann ($\sigma=0$),
and (b) Robin ($\sigma=1$).
 }\label{fig:1D}
 \end{figure}

\begin{figure}[!htbp]
 \centering
  \subfigure{ \includegraphics[scale=.3]{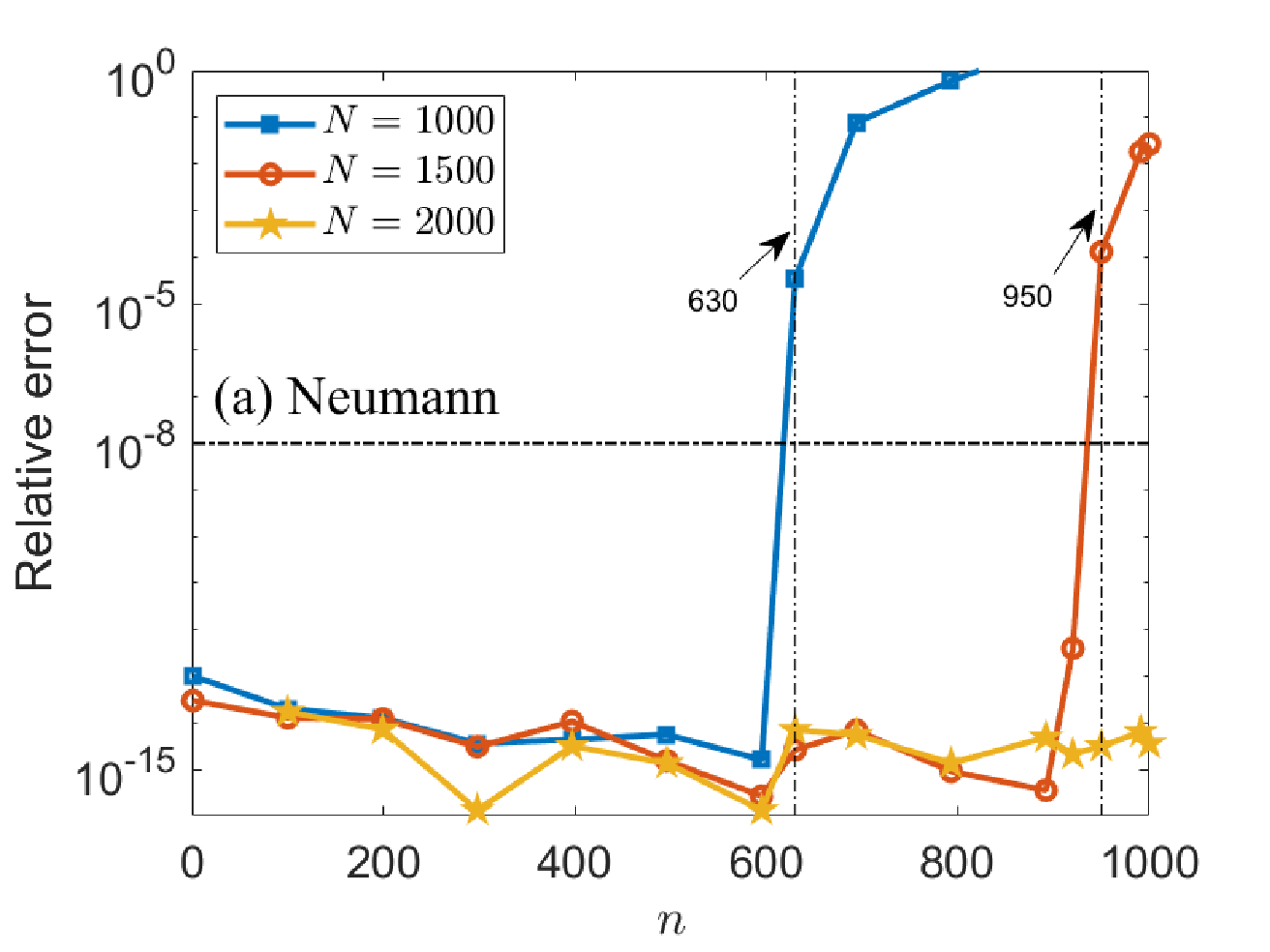}}
  \subfigure{ \includegraphics[scale=.3]{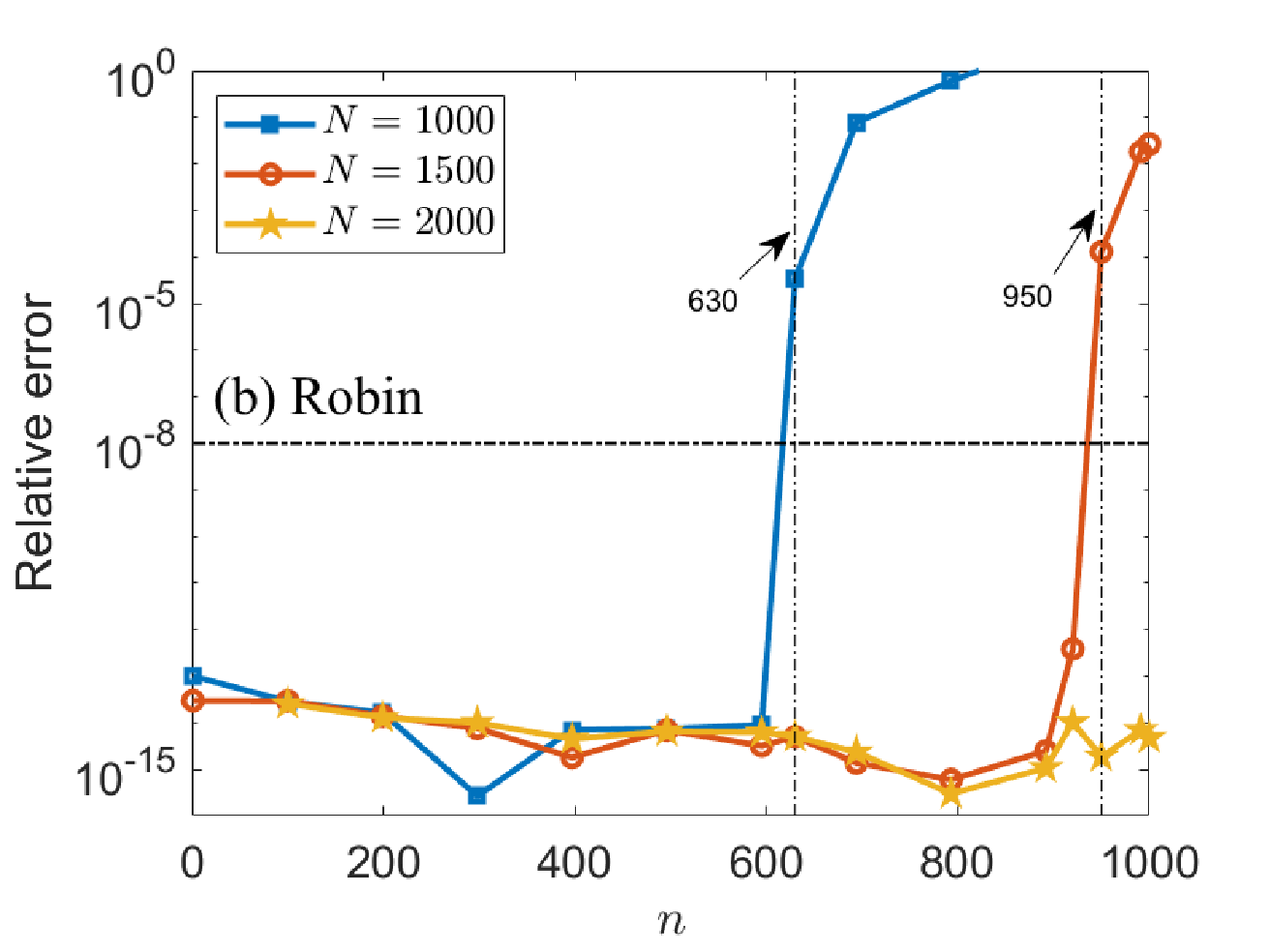}}
  \caption{Accuracy test in 1D with $\Omega=(0,\pi)$ and $V=\frac{x^2}{2}$ in \eqref{SO}-\eqref{SORobin}
  for different $N$ and $\sigma$: (a) Neumann ($\sigma=0$), and (b) Robin ($\sigma=1$).
  }\label{fig:1DV}
 \end{figure}

Figure \ref{fig:1D} displays the relative errors of the eigenvalues
when $V(x)\equiv 0$ for different $\sigma$ and DOFs $N$.
Figure \ref{fig:1DV} plots the relative errors of the eigenvalues
when $V(x)=\frac{x^2}{2}$ for different $\sigma$ and DOFs $N$, while the numerical
reference solutions are obtained by taking $N_{\rm ex}=5000$.

From Figs. \ref{fig:1D}\&\ref{fig:1DV}, we can conclude that,  with a given DOFs $N$, the proposed spectral method
can compute up to approximately $\frac{2}{\pi}N\approx 0.6N$ eigenvalues with relative errors
below $10^{-8}$ (c.f. Figs. \ref{fig:1D}\&\ref{fig:1DV}) for the Schr\"{o}dinger eigenvalue problem \eqref{SO}-\eqref{SORobin} in 1D. With this kind of resolution capacity, we can obtain
numerically the RtN gaps of the Schr\"odinger operator and their convergence rates in 1D, which will
be reported in the next section.

\subsection{Accuracy in 2D}
We take $d=2$ and  $\Omega$ a unit square or a unit disk in
\eqref{SO}-\eqref{SORobin}. When $V(x,y)\equiv 0$,
all the eigenvalues are given analytically under different boundary conditions \cite{gittins2020courant,rudnick2021differences}.
When $\Omega$ is taken as a unit square $(0,1)^2$ and a unit disk, we adopt the
spectral method and the spectral-element method
presented in the previous section to solve this problem numerically for obtaining
the corresponding eigenvalues, respectively.
In the latter case, we also investigate the performance of different combinations
of the number of sub-elements and the polynomial degree of basis functions within
each sub-element, while keeping the total number of degrees of freedom (DOFs) fixed.

\begin{figure}[htbp]
 \centering
  \subfigure{ \includegraphics[scale=.3]{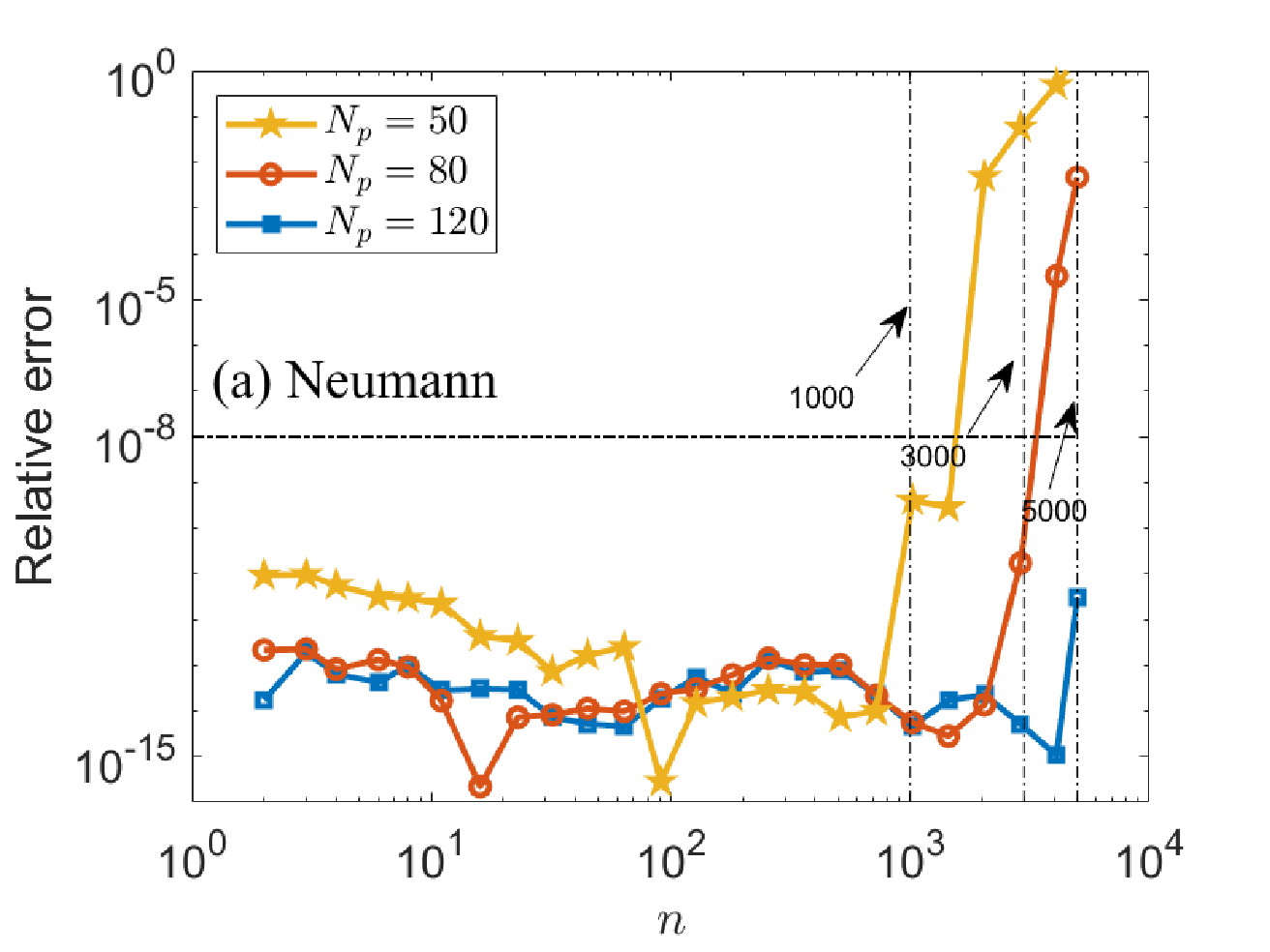}}
  \subfigure{ \includegraphics[scale=.3]{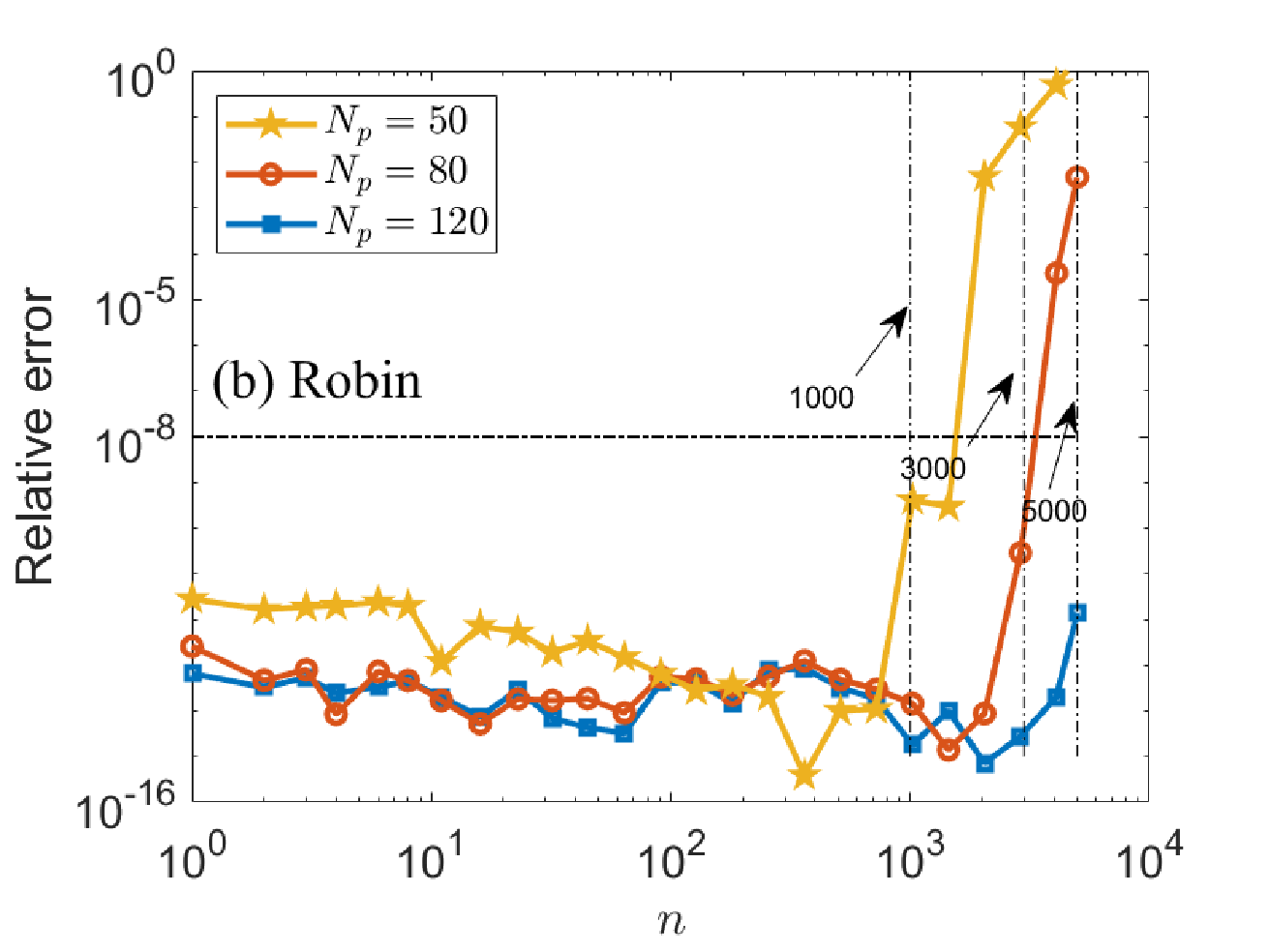}}
  \caption{Accuracy test in 2D with
$\Omega=(0,1)^2$ and  $V(x,y)\equiv0$ in \eqref{SO}-\eqref{SORobin} for different $N=N_p\times N_p$
and $\sigma$: (a) Neumann ($\sigma=0$), and (b) Robin ($\sigma=1$). }\label{fig:Square}
 \end{figure}

\begin{figure}[htbp]
 \centering
  \subfigure{ \includegraphics[scale=.3]{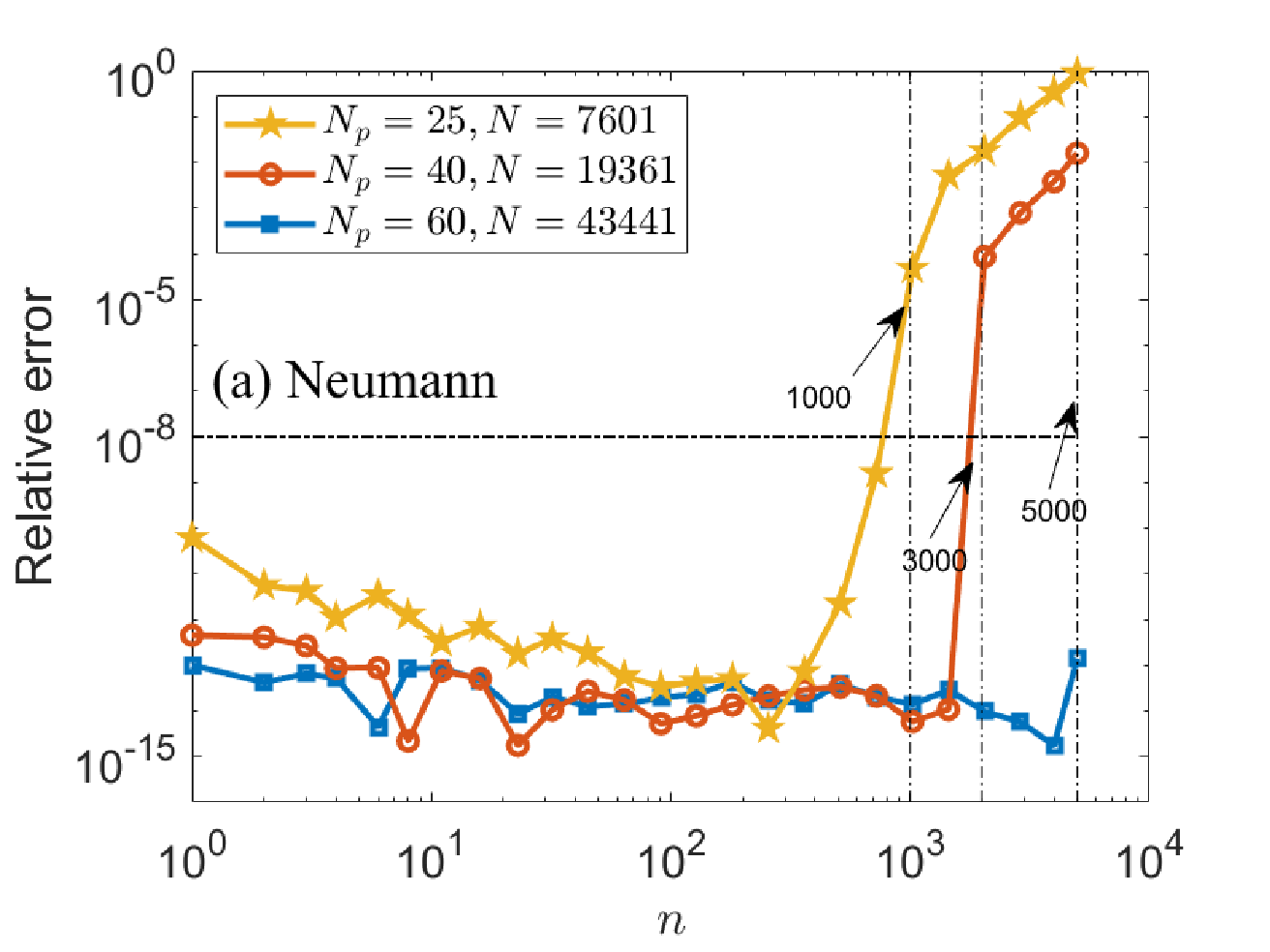}}
  \subfigure{ \includegraphics[scale=.3]{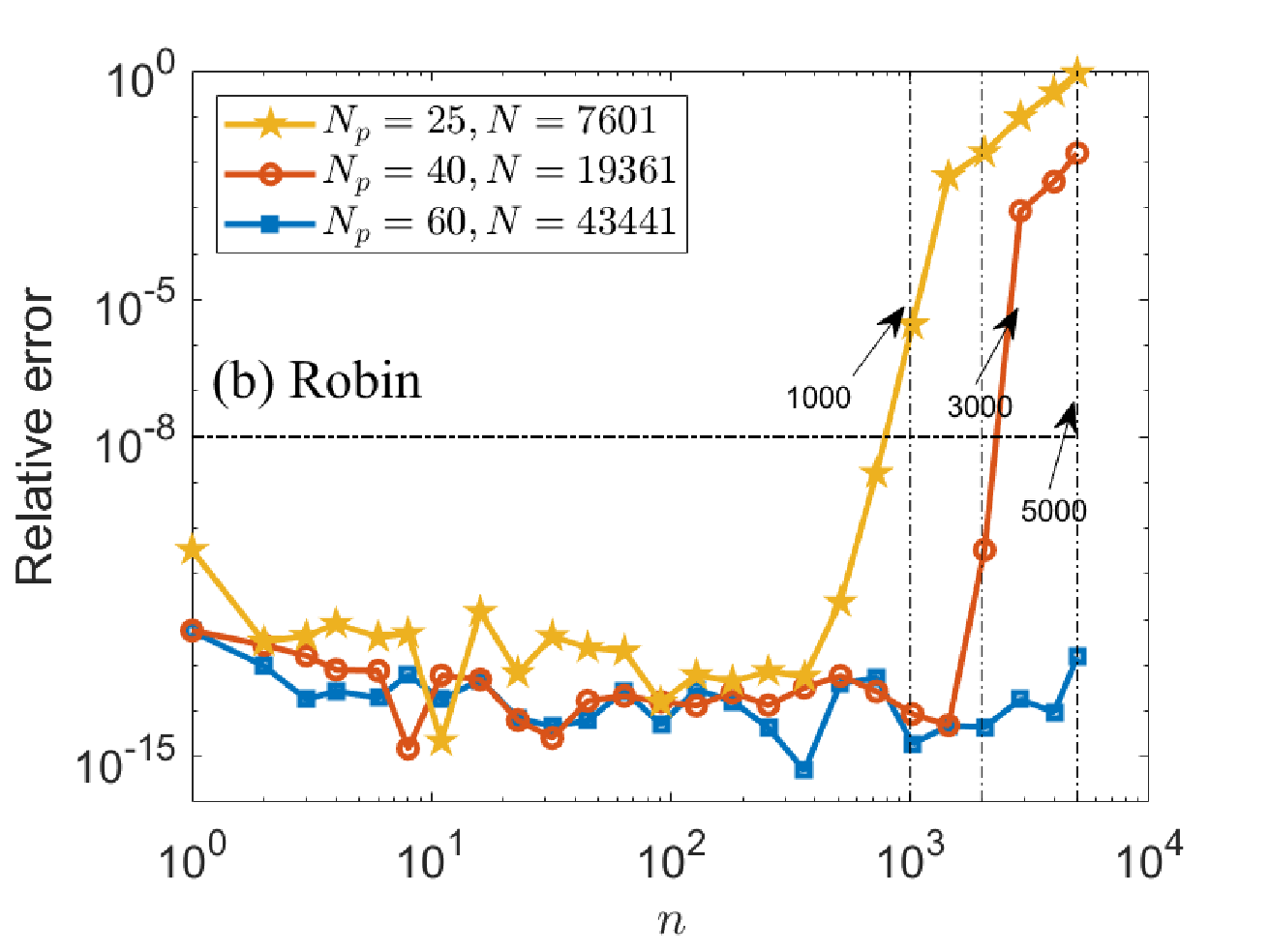}}
  \caption{Accuracy test in 2D with
$\Omega$ taken as a unit disk and  $V(x,y)\equiv0$ in \eqref{SO}-\eqref{SORobin} and $N_E=12$
for different $N_p$ and $\sigma$: (a) Neumann ($\sigma=0$), and (b) Robin ($\sigma=1$). }\label{fig:Disk}
 \end{figure}

\begin{figure}[htbp]
 \centering
  \subfigure{\includegraphics[scale=.4]{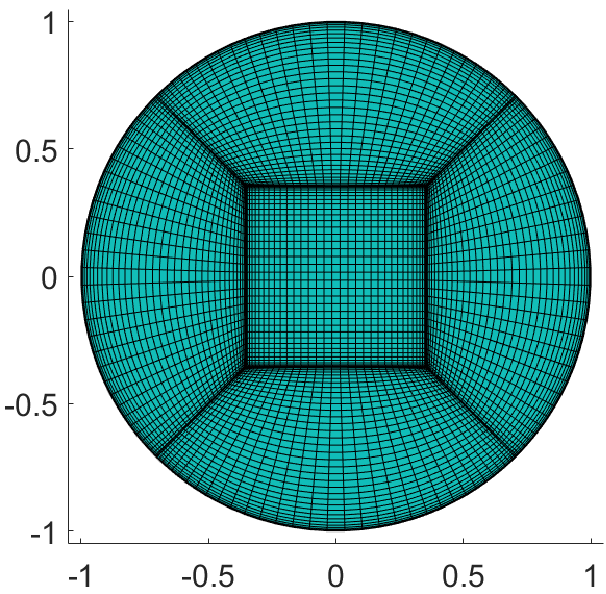}}\qquad \qquad \quad
  \subfigure{\includegraphics[scale=.4]{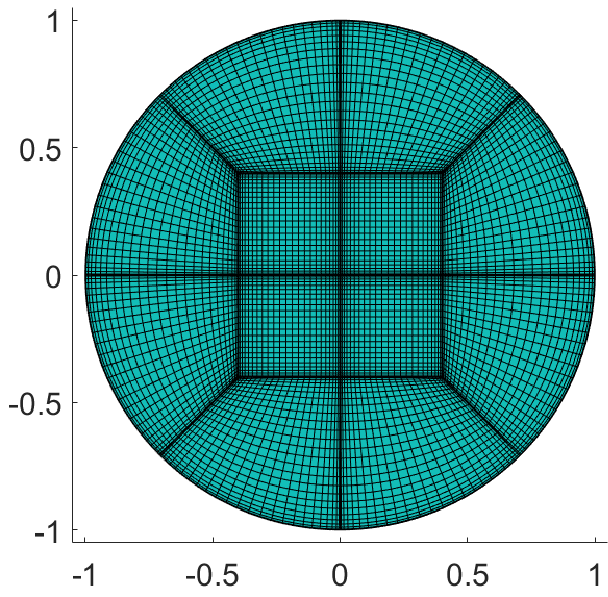}}
  \caption{Schematic illustration of two different spectral-element partitions of the unit disk: $N_E=5$ (left), and
  $N_E=12$ (right).}
  \label{fig:Division}
\end{figure}

\begin{figure}[htbp]
 \centering
  \subfigure{ \includegraphics[scale=.3]{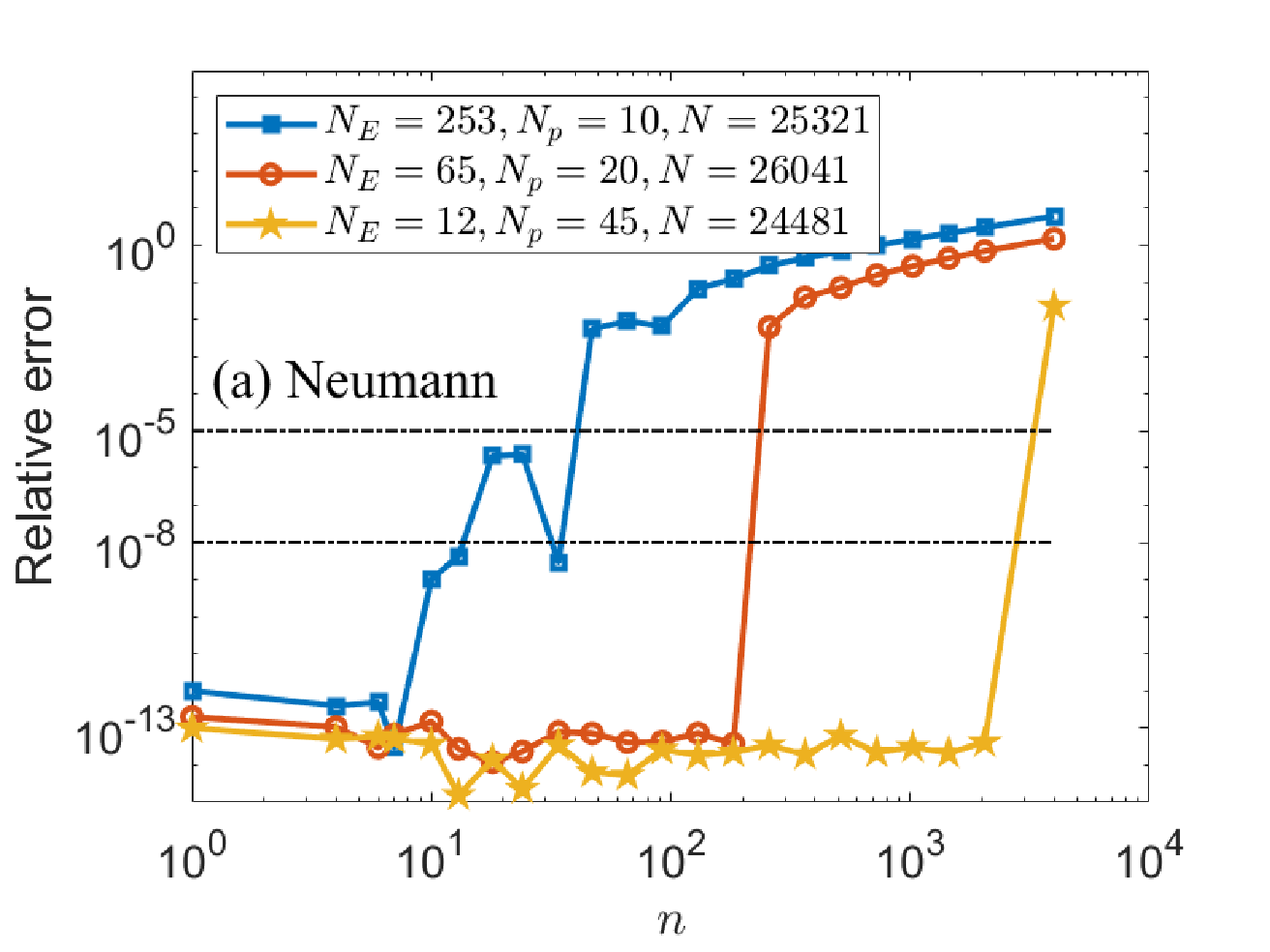}}
  \subfigure{ \includegraphics[scale=.3]{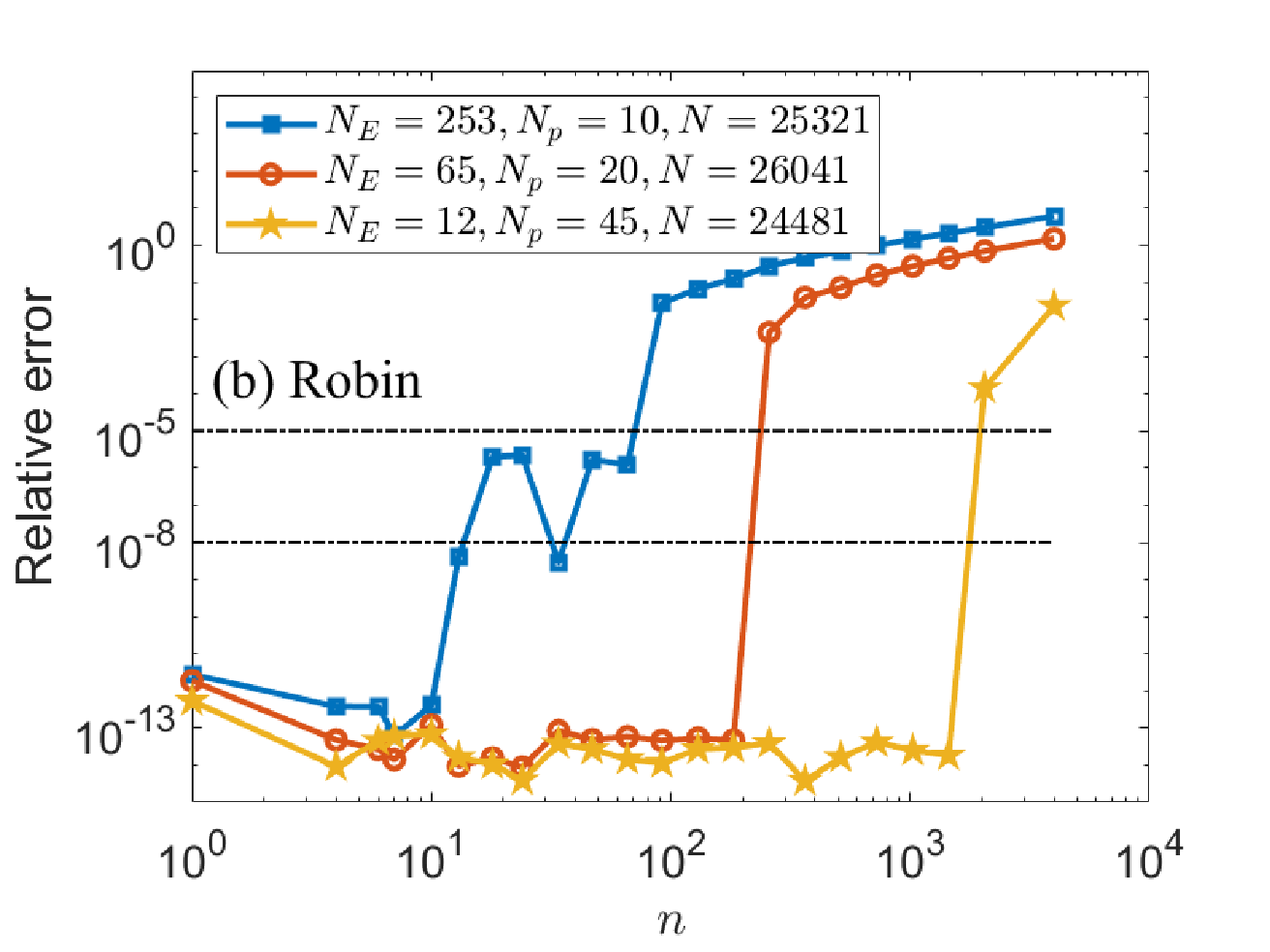}}
  \caption{Accuracy test in 2D with
$\Omega$ taken as a unit disk and  $V(x,y)\equiv0$ in \eqref{SO}-\eqref{SORobin}
for different $N_E$, $N_p$ and $\sigma$: (a) Neumann ($\sigma=0$), and (b) Robin ($\sigma=1$). }\label{fig:hp2}
 \end{figure}

\begin{figure}[htbp]
 \centering
  \subfigure{ \includegraphics[scale=.3]{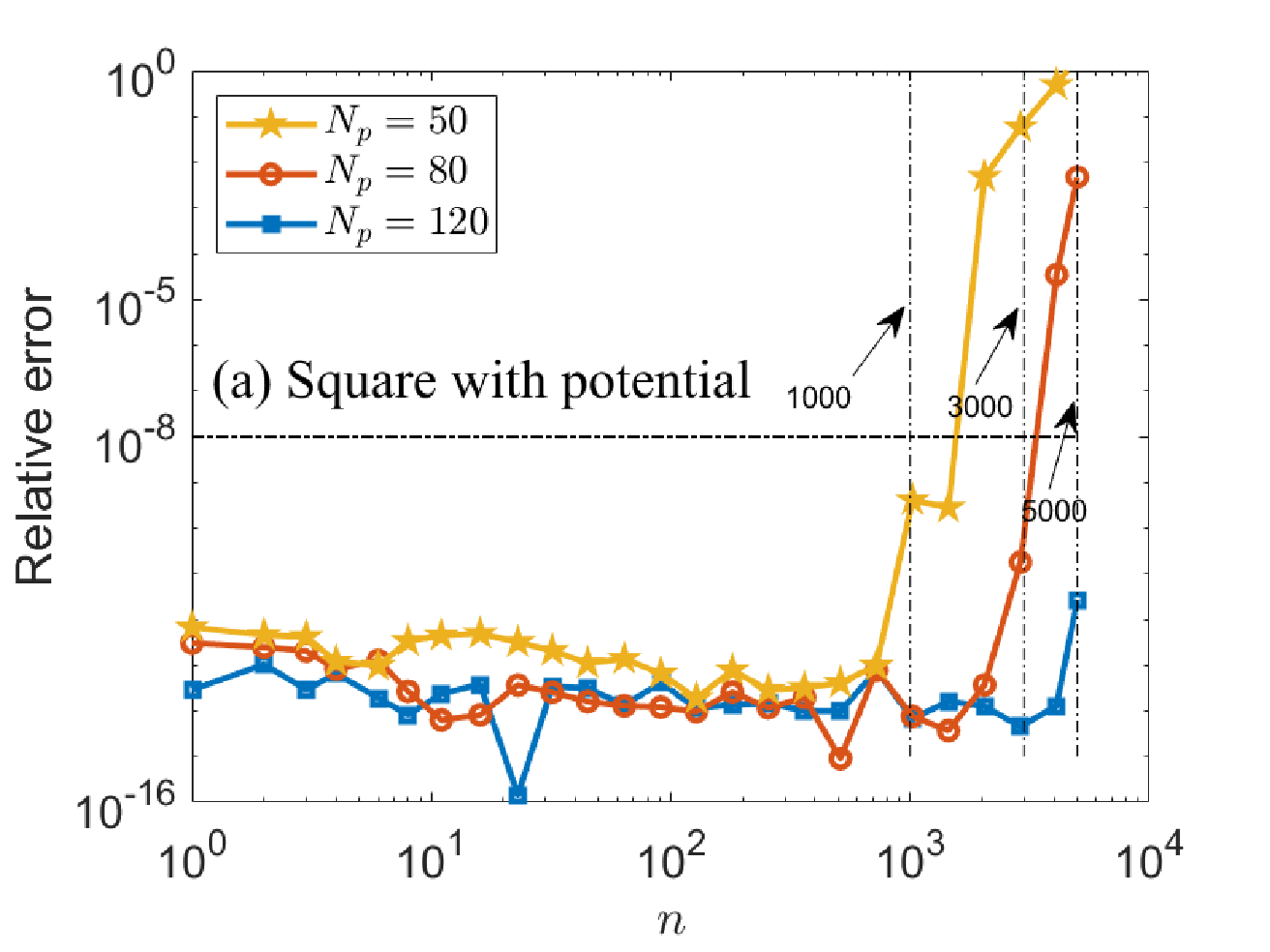}}
  \subfigure{ \includegraphics[scale=.3]{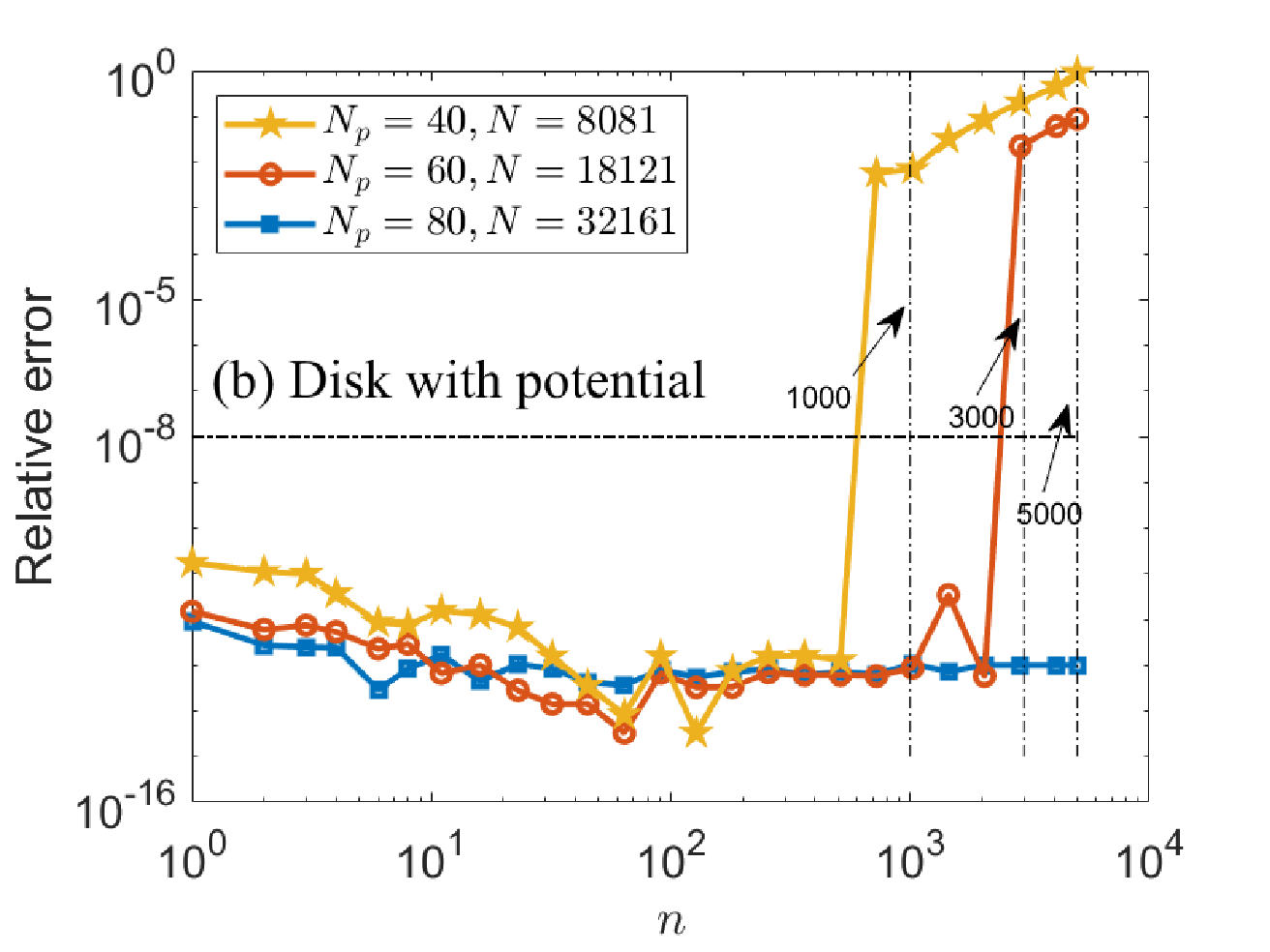}}
  \caption{Accuracy test in 2D with $\sigma=1$ and $V(x,y)=1+\sin(x+y)$ in \eqref{SO}-\eqref{SORobin}
  for different $\Omega$: (a) a square with different $N=N_p\times N_p$, and
 (b) a disk with $N_E=5$ for different $N_p$.}\label{fig:potential}
 \end{figure}

Figure \ref{fig:Square} shows the relative errors of the first $5,000$ eigenvalues
when $V(x,y)\equiv 0$ and $\Omega =(0,1)^2$ in \eqref{SO}-\eqref{SORobin}
for different $\sigma$ and DOFs $N=N_p\times N_p$ with $N_p=50, 80, 120$.
Figure \ref{fig:Disk} plots the relative errors of the first $5,000$ eigenvalues
when $V(x,y)\equiv 0$ and $\Omega$ taken as a unit disk in \eqref{SO}-\eqref{SORobin}
and $N_E=12$
for different $\sigma$ and $N_p$. In addition, Figure \ref{fig:hp2}
displays the relative errors of the eigenvalues
when $V(x,y)\equiv 0$ and $\Omega$ taken as a unit disk in \eqref{SO}-\eqref{SORobin}
for different $\sigma$, $N_E$ and  $N_p$ such that the total DOFs $N\approx 25,000$ (c.f. Fig. \ref{fig:Division}).
Finally, Figure \ref{fig:potential} shows
the relative errors of the first $5,000$ eigenvalues when $\sigma=1$ and $V(x,y)=1+\sin(x+y)$ in \eqref{SO}-\eqref{SORobin} when $\Omega$ is taken as a unit square and a unit disk to be solved
numerically by the spectral method with different $N=N_P\times N_p$  and
the spectral-element method with different $N_E$ and $N_P$, respectively.

From Figs. \ref{fig:Square}-\ref{fig:potential}, we can draw the following conclusions:
(i) With the same number of DOFs, we observe that a larger number of reliable eigenvalues
can generally be obtained by increasing the polynomial degree $N_P$ and decreasing
the number of elements $N_E$ (c.f. Fig. \ref{fig:hp2}), which
demonstrates the significant advantage of high-order spectral methods in computing a large number of eigenvalues with high accuracy. (ii) With a given DOFs $N$, the proposed spectral method
can compute up to approximately $\frac{4}{\pi^2}N\approx 0.4N$ eigenvalues with relative errors
below $10^{-8}$ (c.f. Figs. \ref{fig:Square}\&\ref{fig:potential} (left one)) for the Schr\"{o}dinger eigenvalue problem \eqref{SO}-\eqref{SORobin} in 2D. (iii) The proposed spectral-element method
can compute up to more than $5,000$ eigenvalues with relative errors
below $10^{-8}$ (c.f. Figs. \ref{fig:Disk}\&\ref{fig:potential} (right one)) for the Schr\"{o}dinger eigenvalue problem \eqref{SO}-\eqref{SORobin} in 2D. (iv) In 2D, with about DOFs at $N\approx 25,000$, we can obtain
numerically at least $5,000$ reliable eigenvalues (with relative error less than $10^{-8}$) (c.f. Figs. \ref{fig:Square}, \ref{fig:Disk} and \ref{fig:potential}). Again, with
this kind of resolution capacity of the proposed spectral and spectral-element methods, we can obtain
numerically the RtN gaps of the Schr\"odinger operator and their convergence rates in 2D, which will
be reported in the next section.

\begin{figure}[htbp]
 \centering
  \subfigure{ \includegraphics[scale=.3]{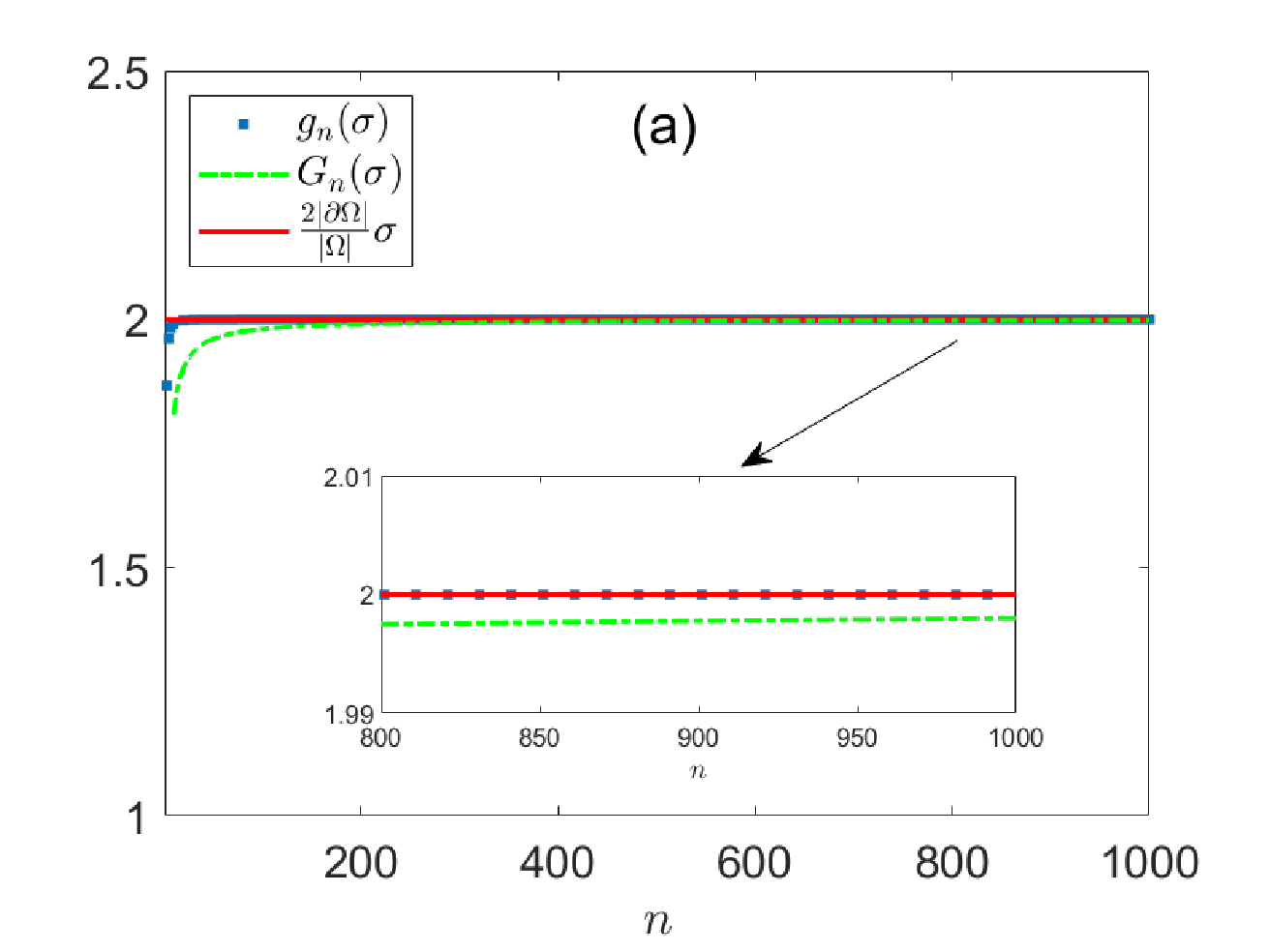}}\qquad
  \subfigure{ \includegraphics[scale=.3]{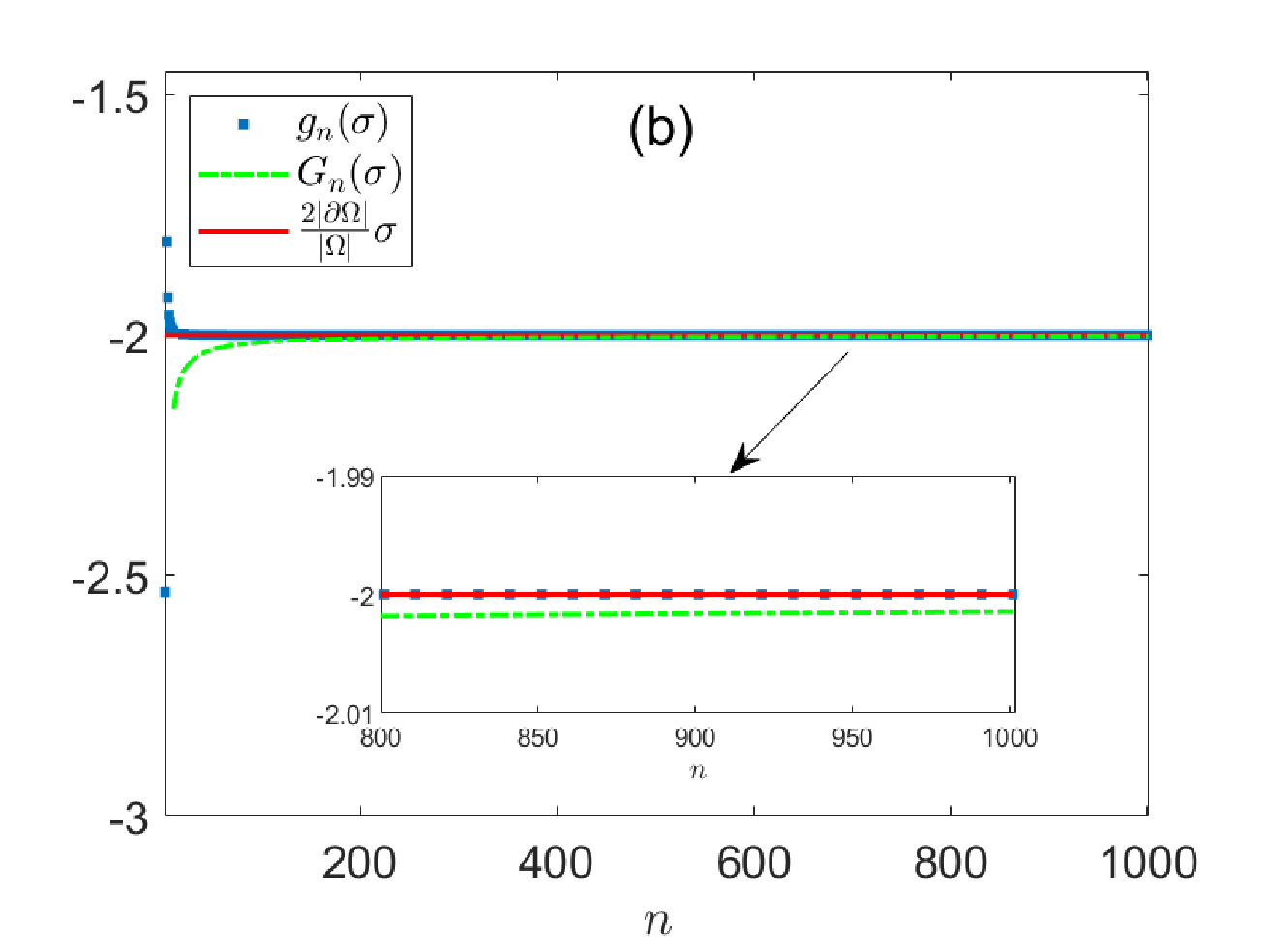}}
  \caption{The RtN gaps of \eqref{SO}-\eqref{SORobin} in 1D with $V=\frac{x^2}{2}$ for different $\sigma$:
  (a) $\sigma=1$, and (b) $\sigma=-1$. }\label{fig:RNgap1D}
 \end{figure}

\section{Numerical studies on the Robin-to-Neumann (RtN) gaps}
In this section, we investigate numerically the Robin-to-Neumann (RtN) gaps in 1D, 2D and 3D by
using the spectral and spectral-element methods proposed in Section 3. In each case, we take enough large
DOFs $N$ such that we have enough reliable eigenvalues for obtaining conclusions on the RtN gaps.

\subsection{The RtN gaps in 1D}
We take $d=1$ and  $\Omega=(-1,1)$ in \eqref{SO}-\eqref{SORobin}. The problem
is solved numerically by the spectral method.
Figure \ref{fig:RNgap1D} shows the RtN gaps with $V(x)=\frac{x^2}{2}$ for different $\sigma$.

 \begin{figure}[htbp]
 \centering
  \subfigure{ \includegraphics[scale=.35]{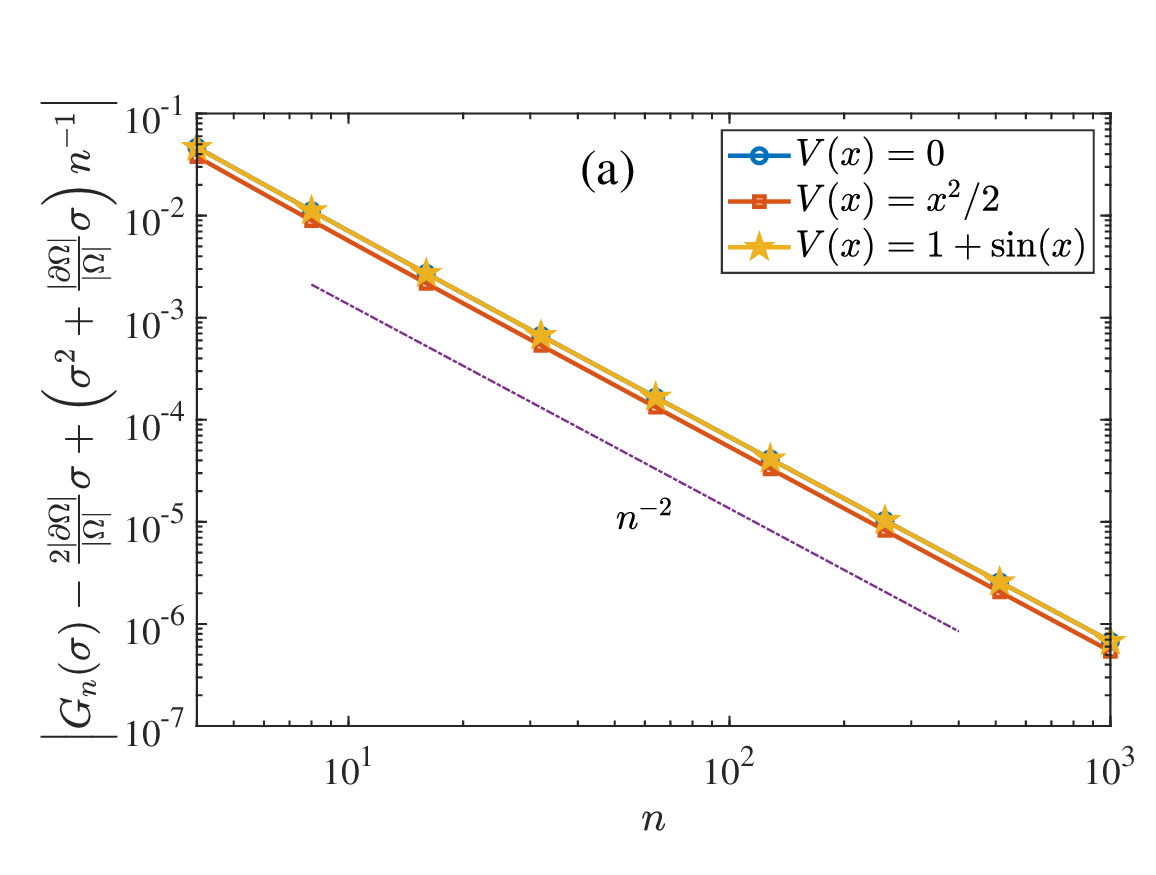}}\qquad
  \subfigure{ \includegraphics[scale=.35]{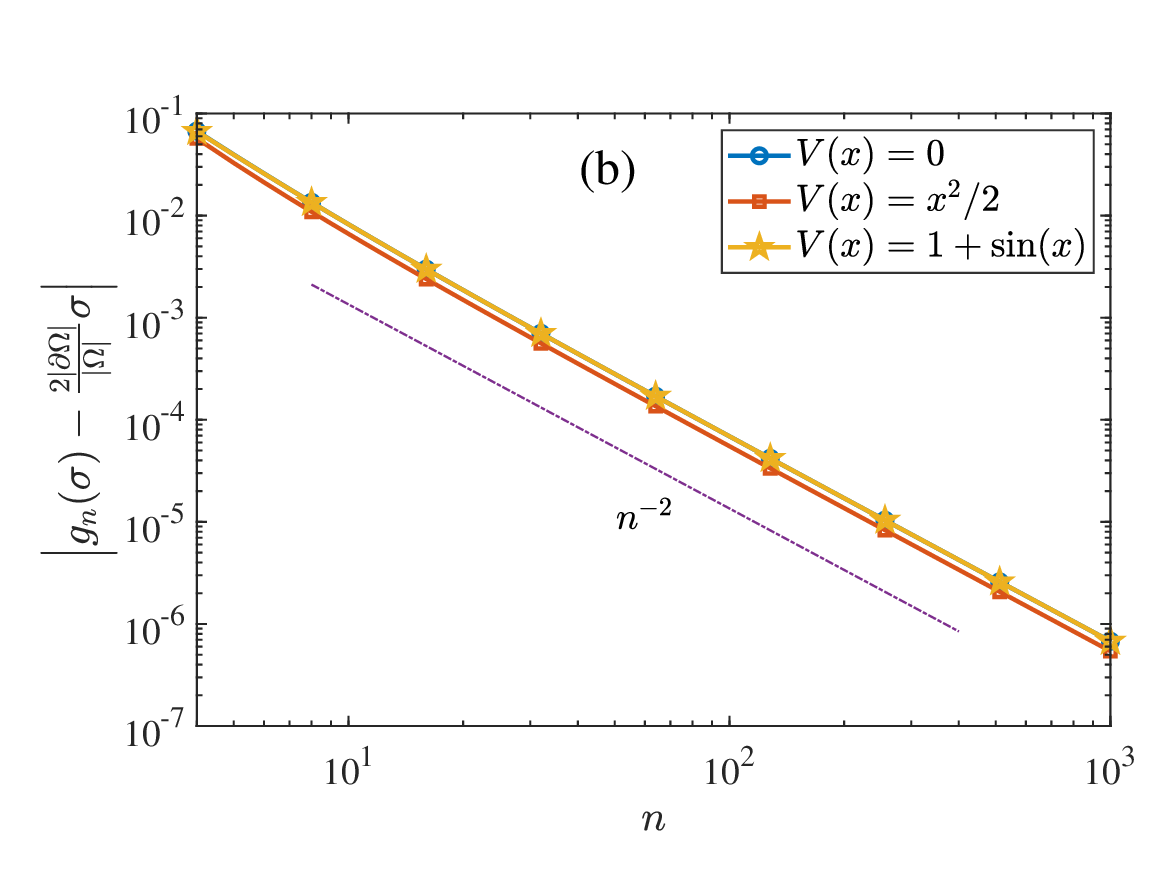}}
  \caption{Convergence rates of the RtN gaps of \eqref{SO}-\eqref{SORobin} in 1D with $\sigma=1$
  for different potentials $V(x)$: (a) convergence rate of the cumulative averages $G_n(\sigma)$, and
  (b) convergence rate of the RtN gaps $g_n(\sigma)$. }\label{fig:Conver1D}
 \end{figure}

\begin{figure}[!htbp]
 \centering
  \subfigure{ \includegraphics[scale=.35]{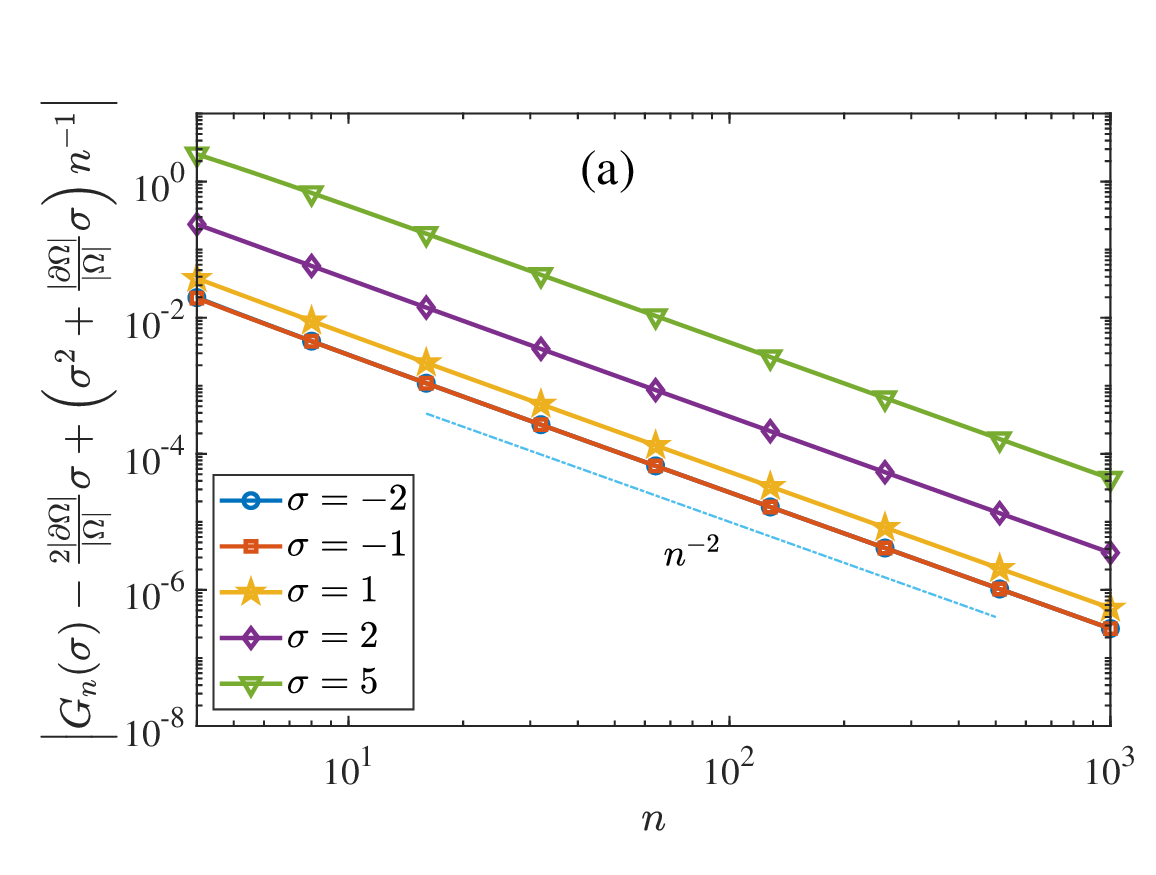}}\qquad
  \subfigure{ \includegraphics[scale=.35]{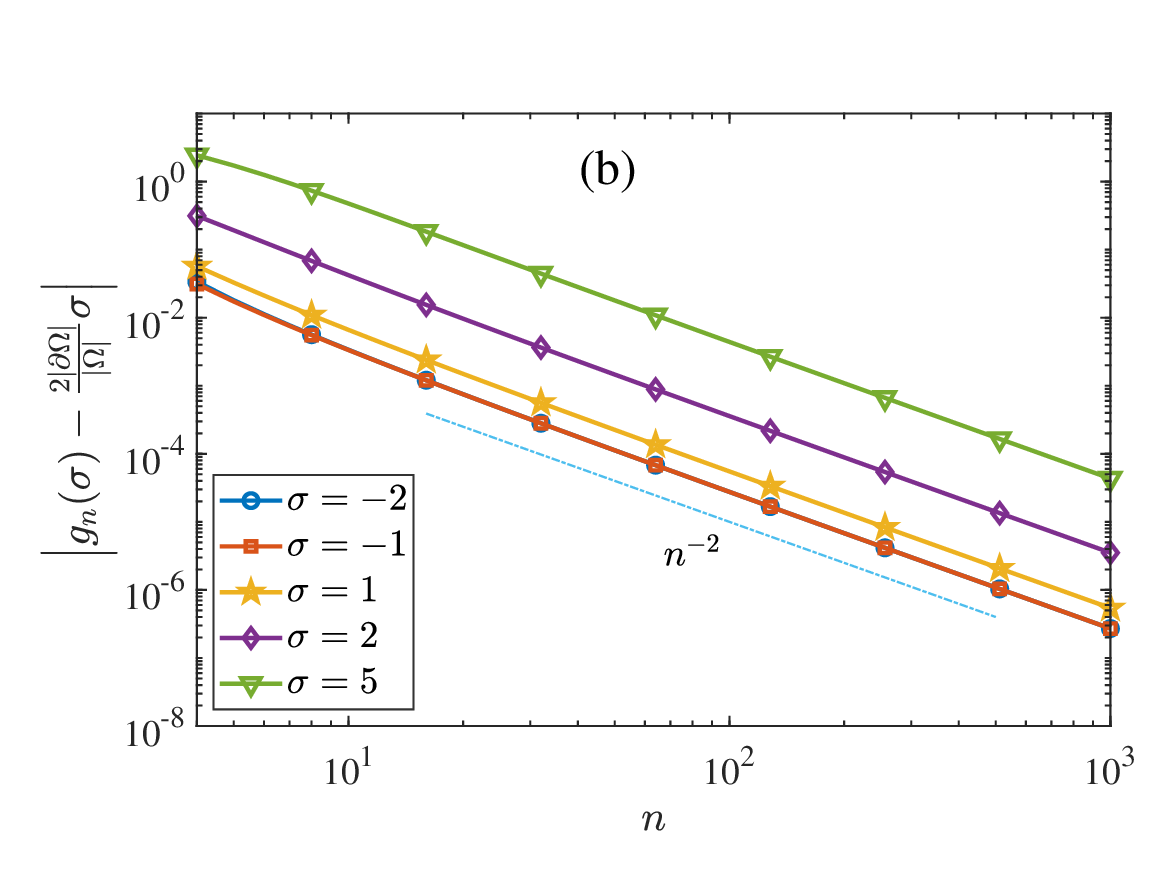}}
  \caption{Convergence rates of the RtN gaps of \eqref{SO}-\eqref{SORobin} in 1D with $V(x)=\frac{x^2}{2}$
  for different $\sigma$: (a) convergence rate of the cumulative averages $G_n(\sigma)$, and
  (b) convergence rate of the RtN gaps $g_n(\sigma)$. }\label{fig:diff_sigma1D}
 \end{figure}

From Figs. \ref{fig:RNgap1D}-\ref{fig:diff_sigma1D}, and additional results
with several other potentials, such as $V(x)=1+\sin(x)$, or $V(x)=1-\cos(x)$, or
\[V(x)=\begin{cases} 1, & -\frac{1}{2} \le x \le \frac{1}{3},\\ 0,
& {\rm otherwise}, \end{cases}, \quad -1 \le x \le 1;
\]
which are not shown here for brevity,
we can immediately conclude \eqref{RtNgapdim} when $d=1$
and \eqref{RtNgap1d}.

\begin{figure}[htbp]
 \centering
 \subfigure{ \includegraphics[scale=.3]{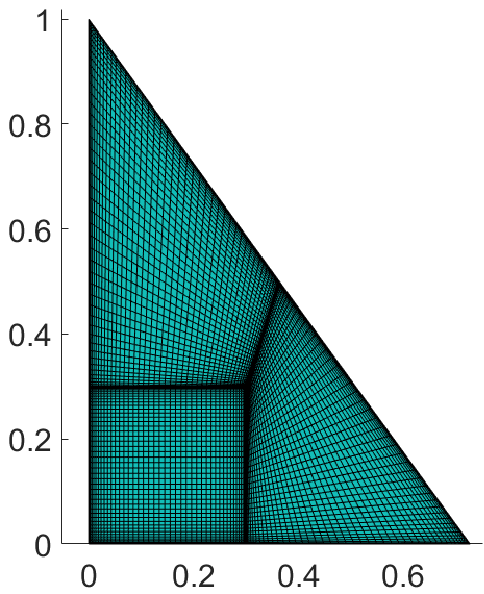}}\quad
 \subfigure{ \includegraphics[scale=.3]{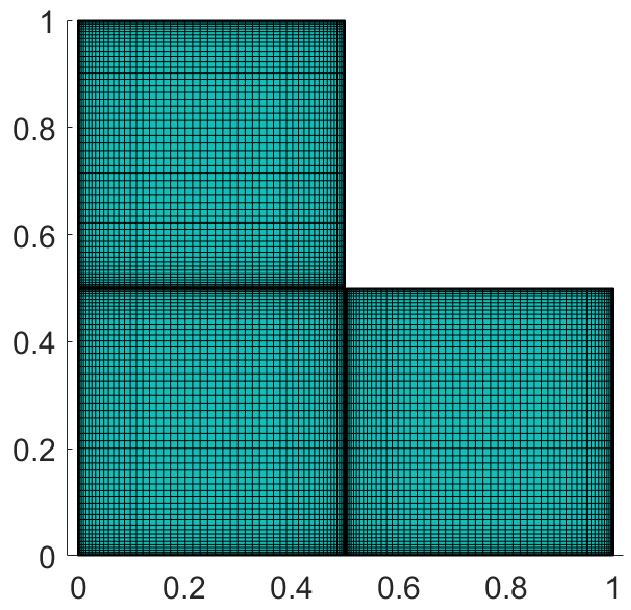}}\quad
  \subfigure{ \includegraphics[scale=.3]{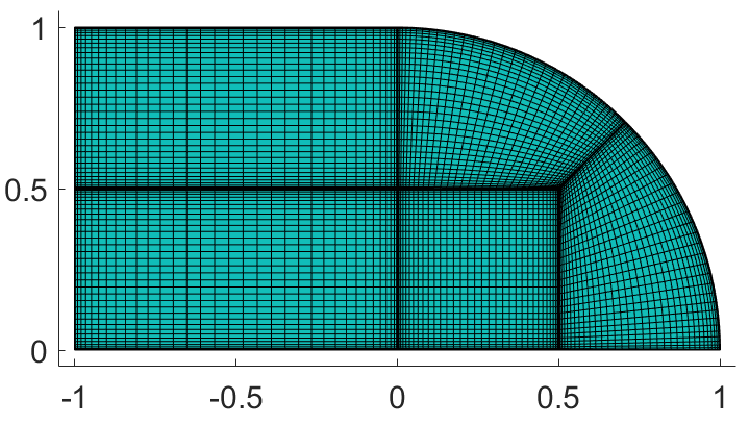}}\quad
  \subfigure{ \includegraphics[scale=.3]{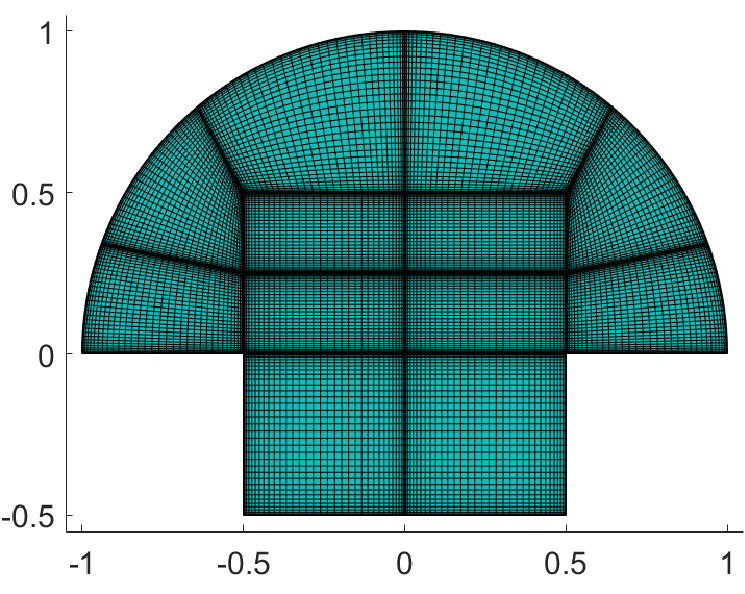}}\quad
  \caption{ Schematics and partitioning of the different domains $\Omega$ in 2D:  a right triangle with an angle $\pi/5$, an L-shaped domain, a quarter stadium, a mushroom (from left to right).}\label{fig:division2}
 \end{figure}
\begin{figure}[htbp]
 \centering
  \subfigure{ \includegraphics[scale=.32]{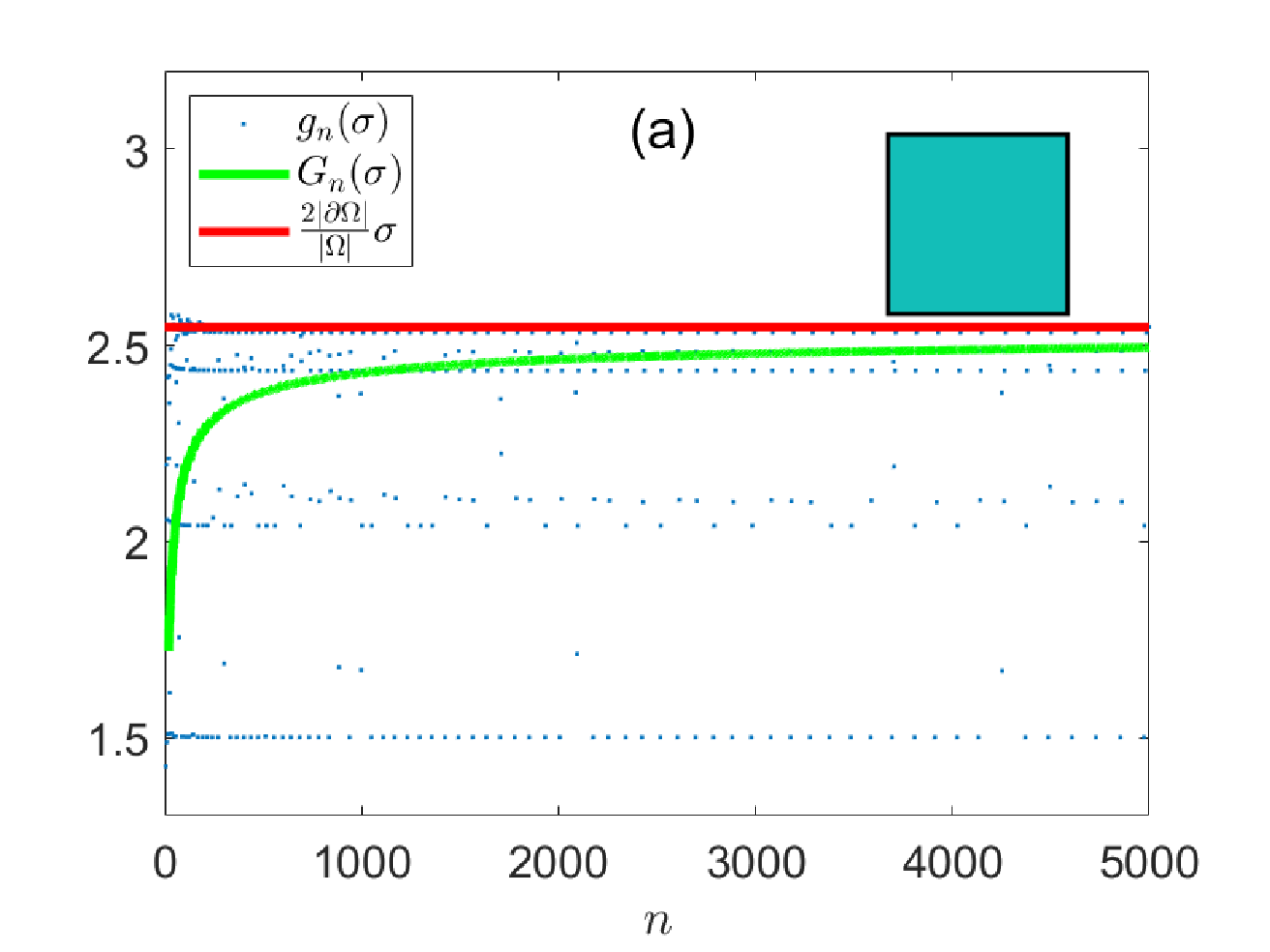}}
  \subfigure{ \includegraphics[scale=.32]{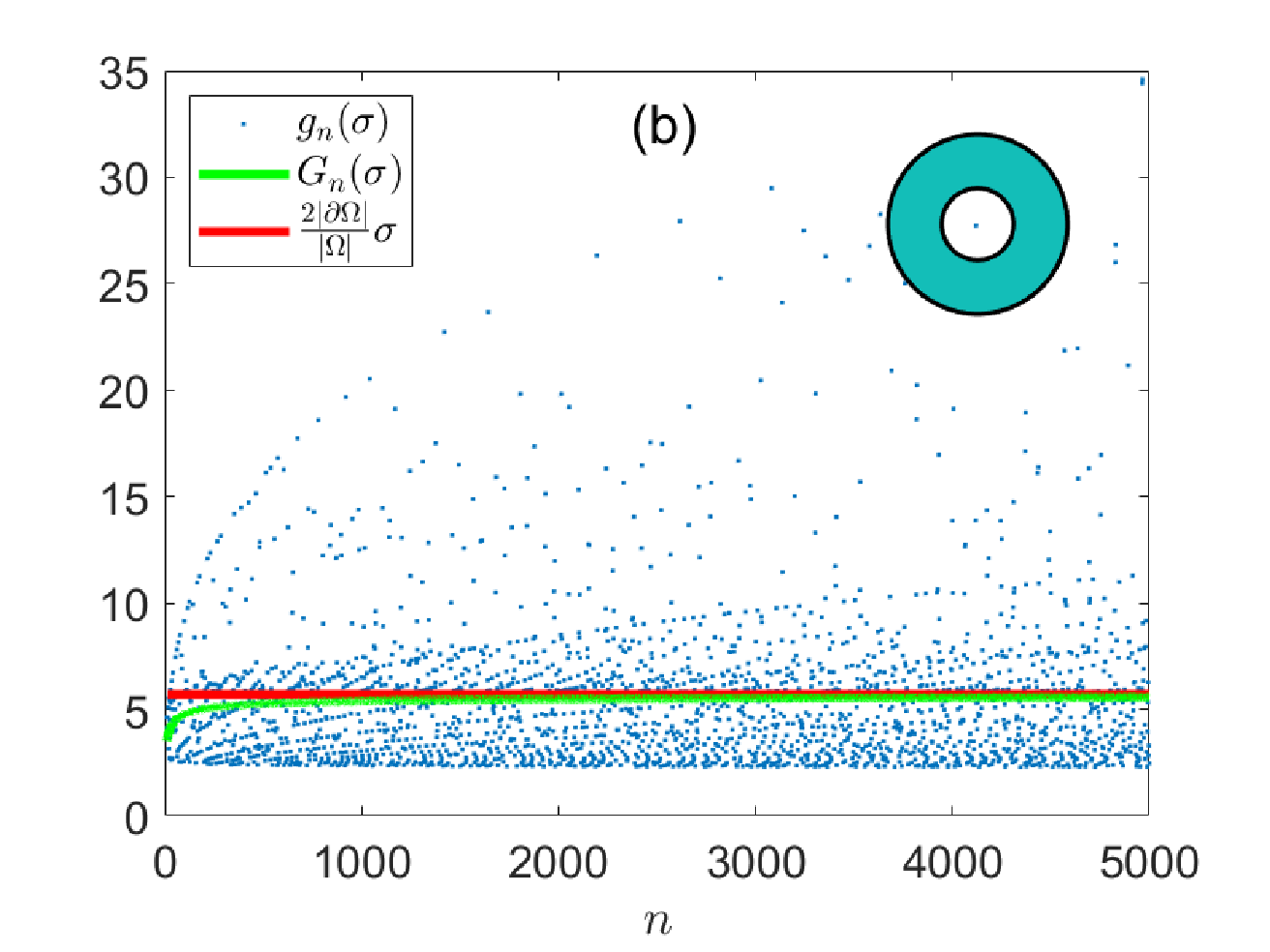}}\\
  \subfigure{ \includegraphics[scale=.32]{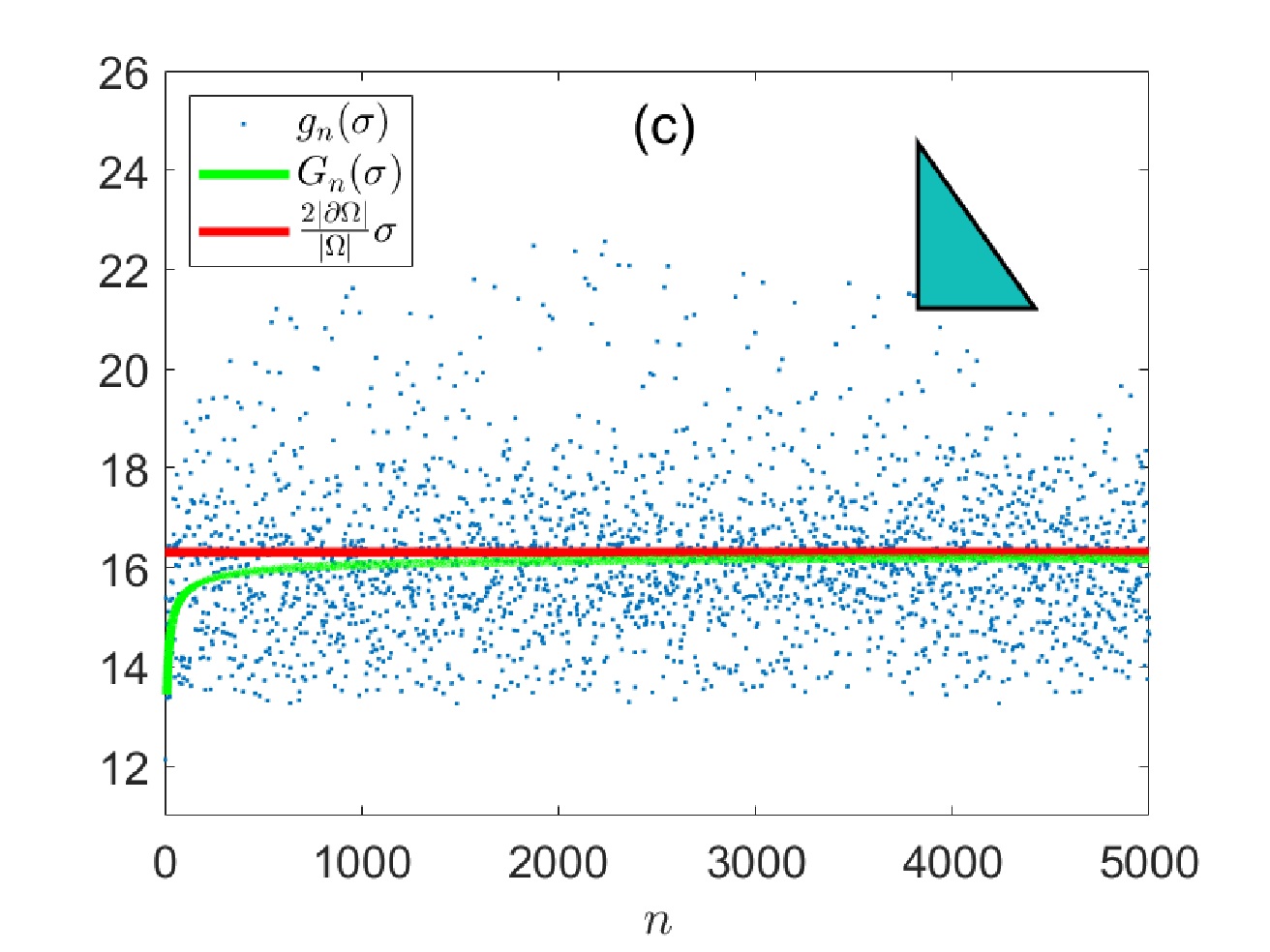}}
  \subfigure{ \includegraphics[scale=.32]{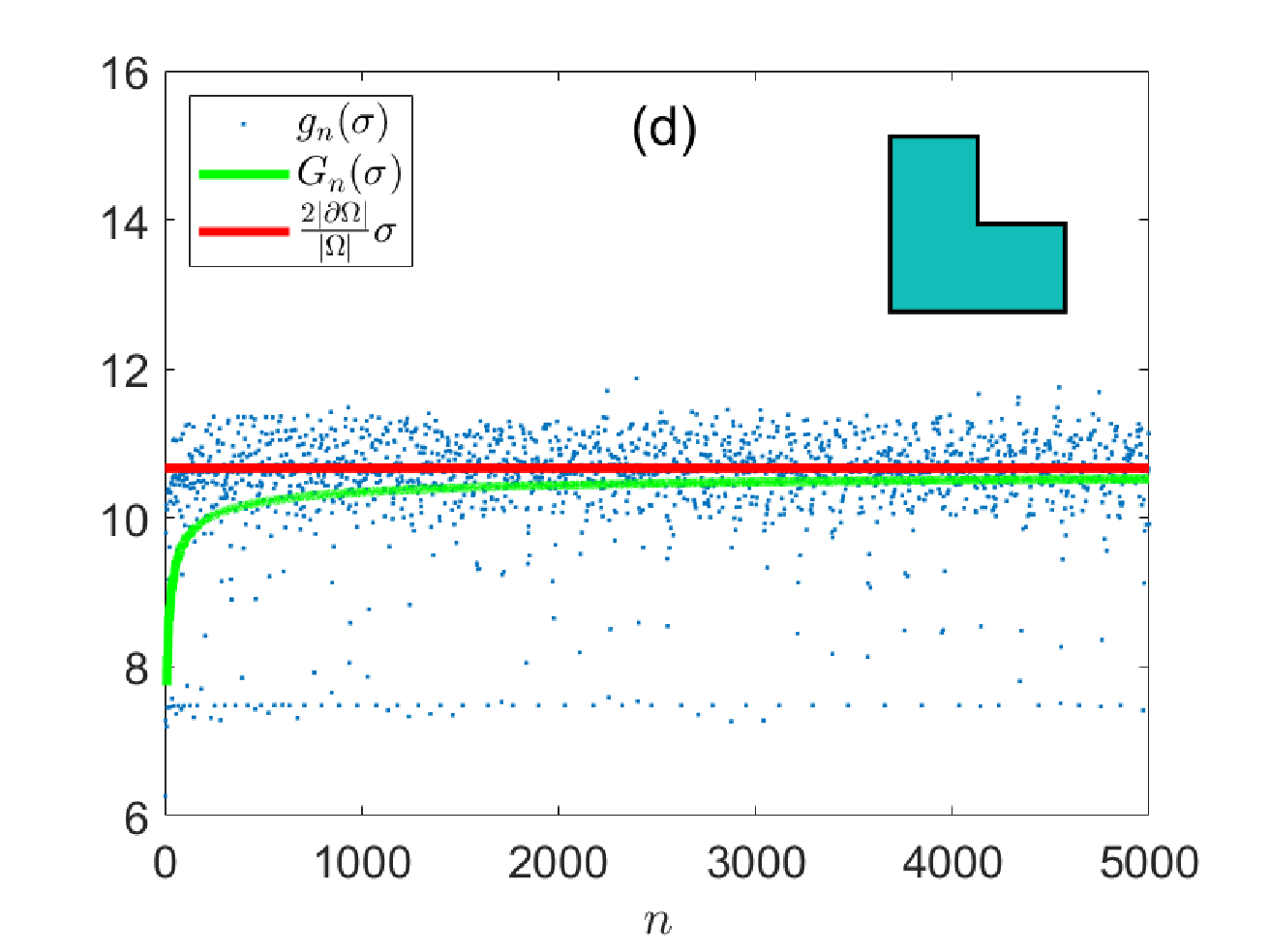}}\\
    \subfigure{ \includegraphics[scale=.32]{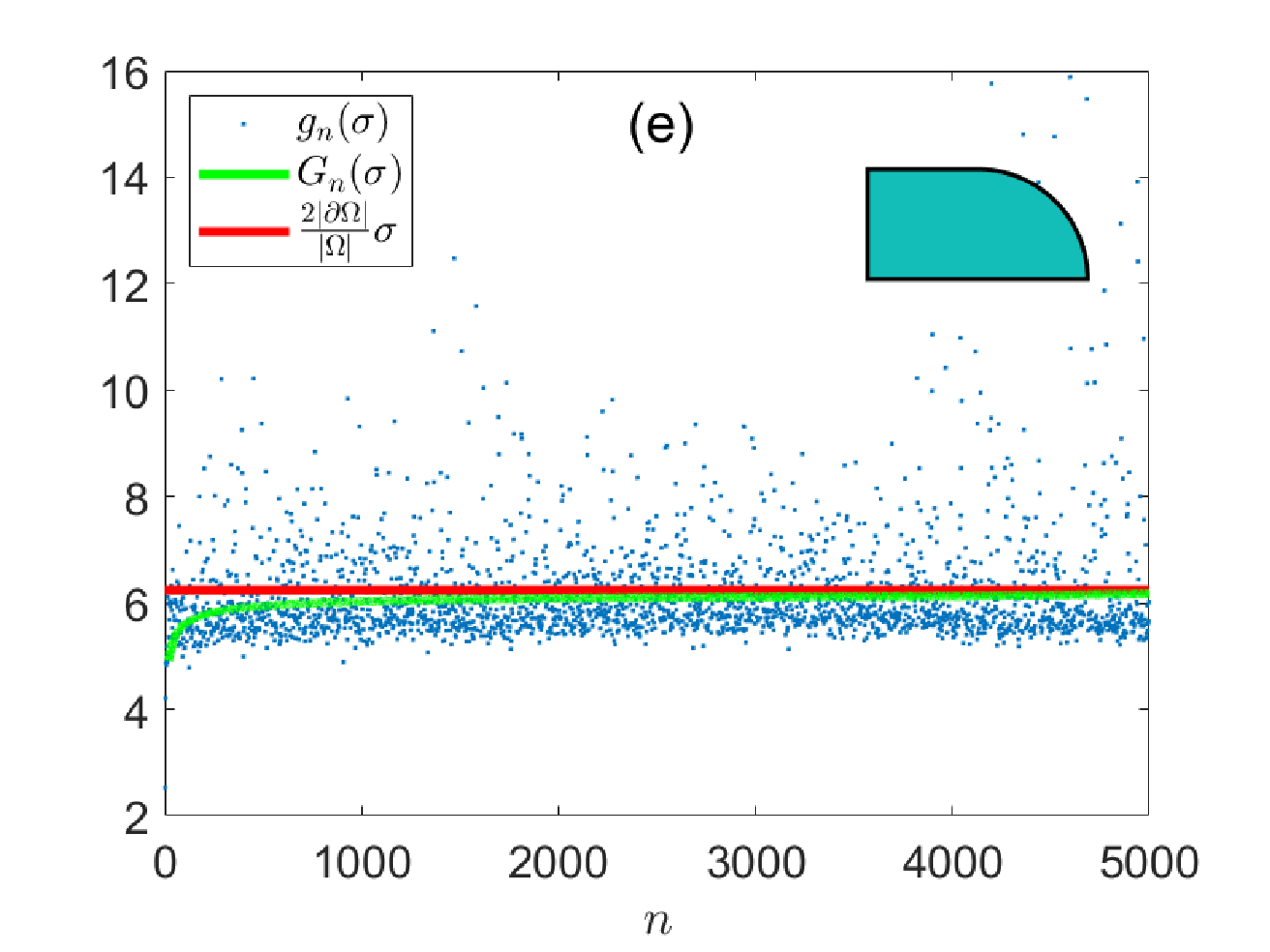}}
      \subfigure{ \includegraphics[scale=.32]{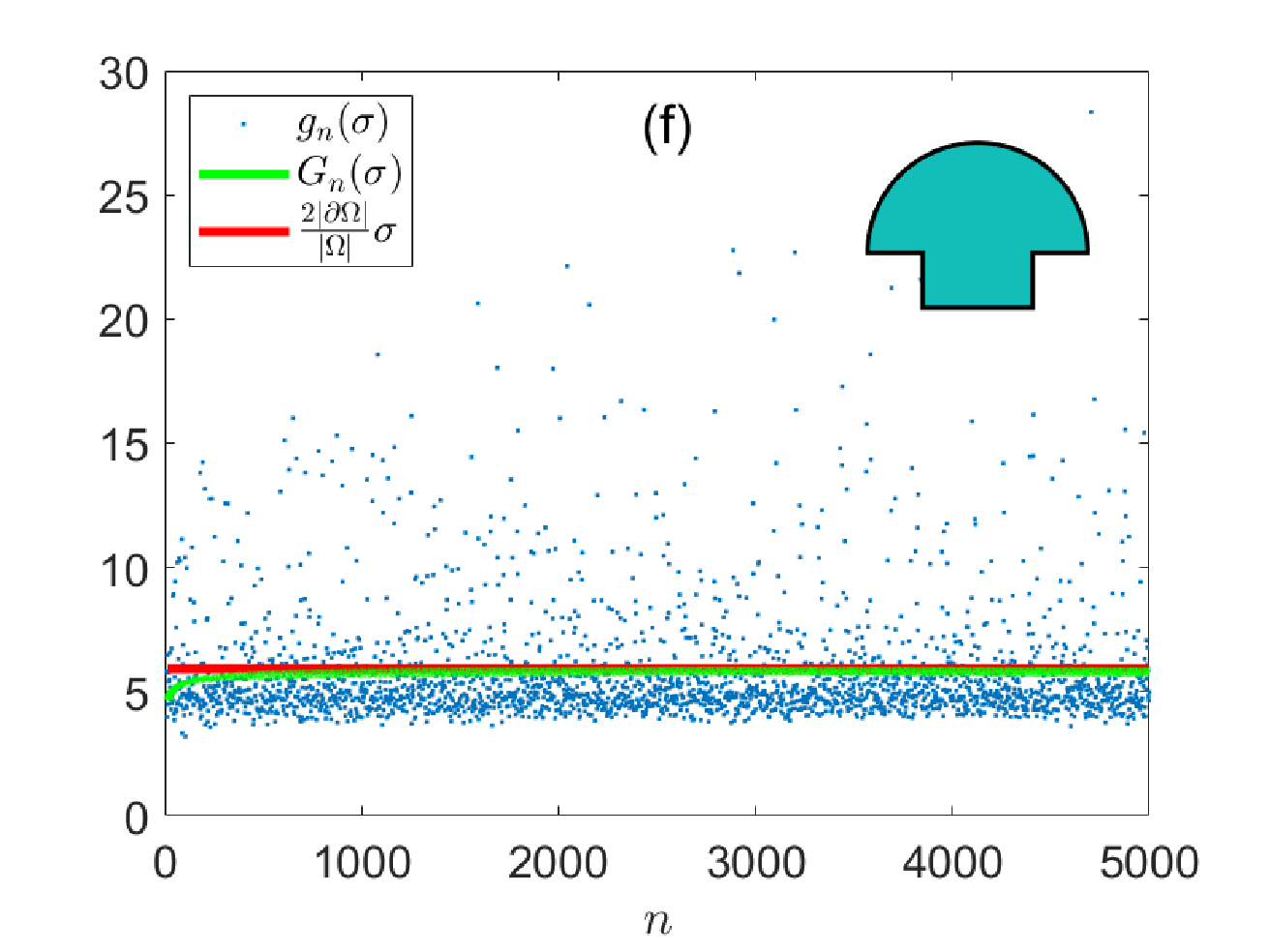}}
  \caption{The RtN gaps of \eqref{SO}-\eqref{SORobin} in 2D with
  $V(x,y)=x^2+y^2$ and $\sigma=1$ for different domains $\Omega$:
  (a) the unit square, (b) the annulus with inner and outer radius as $1$ and $2$, respectively,
  (c) the right triangle with an angle $\frac{\pi}{5}$, (d) the L-shaped domain, (e) the quarter stadium, and  (f) the mushroom.}\label{fig:RNgapV1}
 \end{figure}

\subsection{The RtN gaps in 2D}
We take $d=2$ in \eqref{SO}-\eqref{SORobin} with several different domains, such as
right triangle with an angle $\pi/5$, L-shaped domain, quarter stadium, mushroom (c.f. Fig. \ref{fig:division2}).
The spectral-element method is adopted with $N_E$, $N_p$ and DOFs $N$ given below (c.f. Fig. \ref{fig:division2}):
 \begin{itemize}
 \item Triangle: $N_E=3$, $N_p=128$ and the DOF $N=49,537$;
 \item L-shaped domain: $N_E=3$, $N_p=128$ and the DOF $N=49,665$;
 \item Quarter stadium: $N_E=5$, $N_p=100$ and the DOF $N=50,401$;
 \item Mushroom: $N_E=12$, $N_p=70$ and the DOF $N=59,221$.
 \end{itemize}

\begin{figure}[htbp]
 \centering
  \subfigure{ \includegraphics[scale=.32]{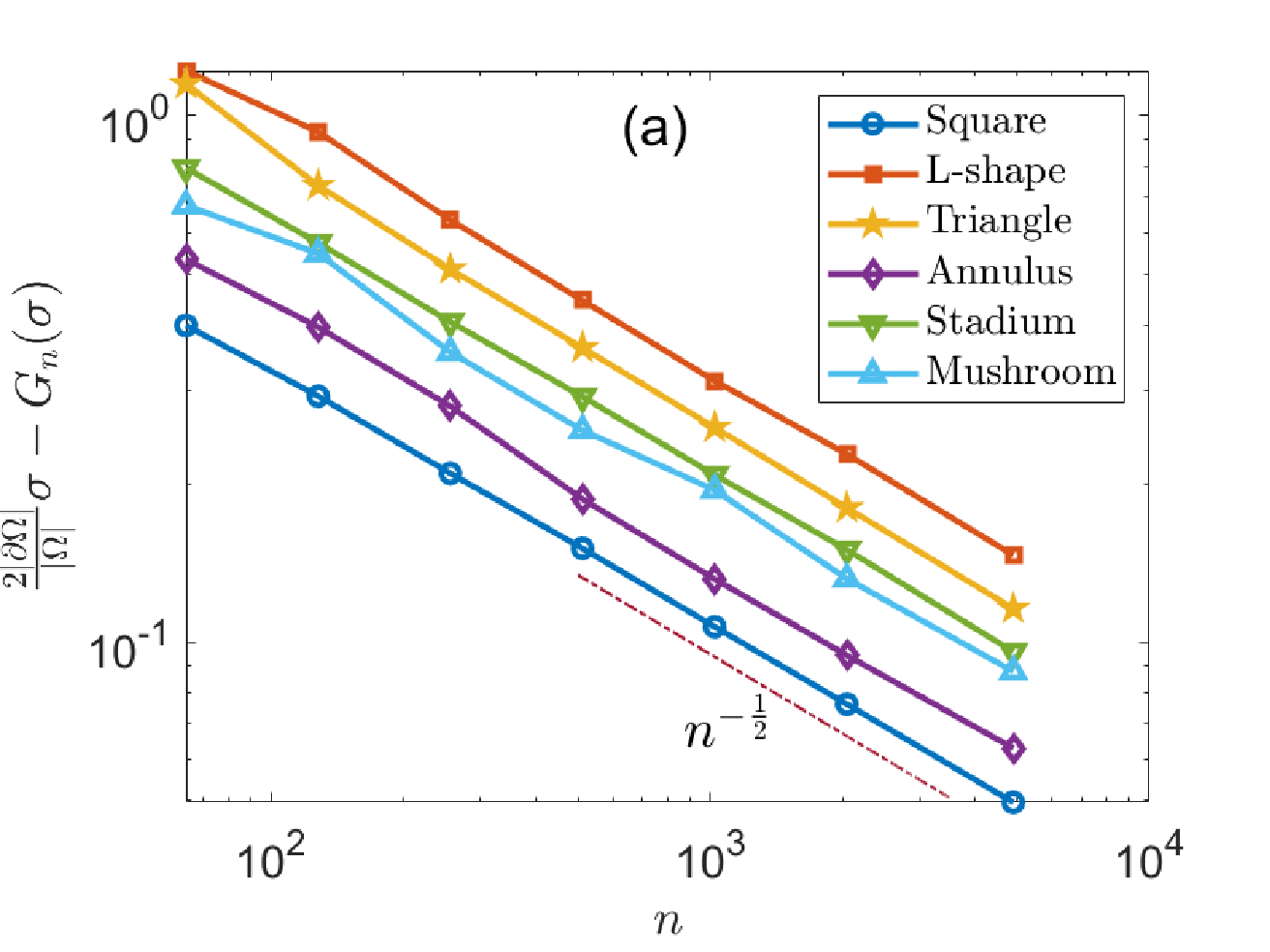}}\hfill
  \subfigure{ \includegraphics[scale=.32]{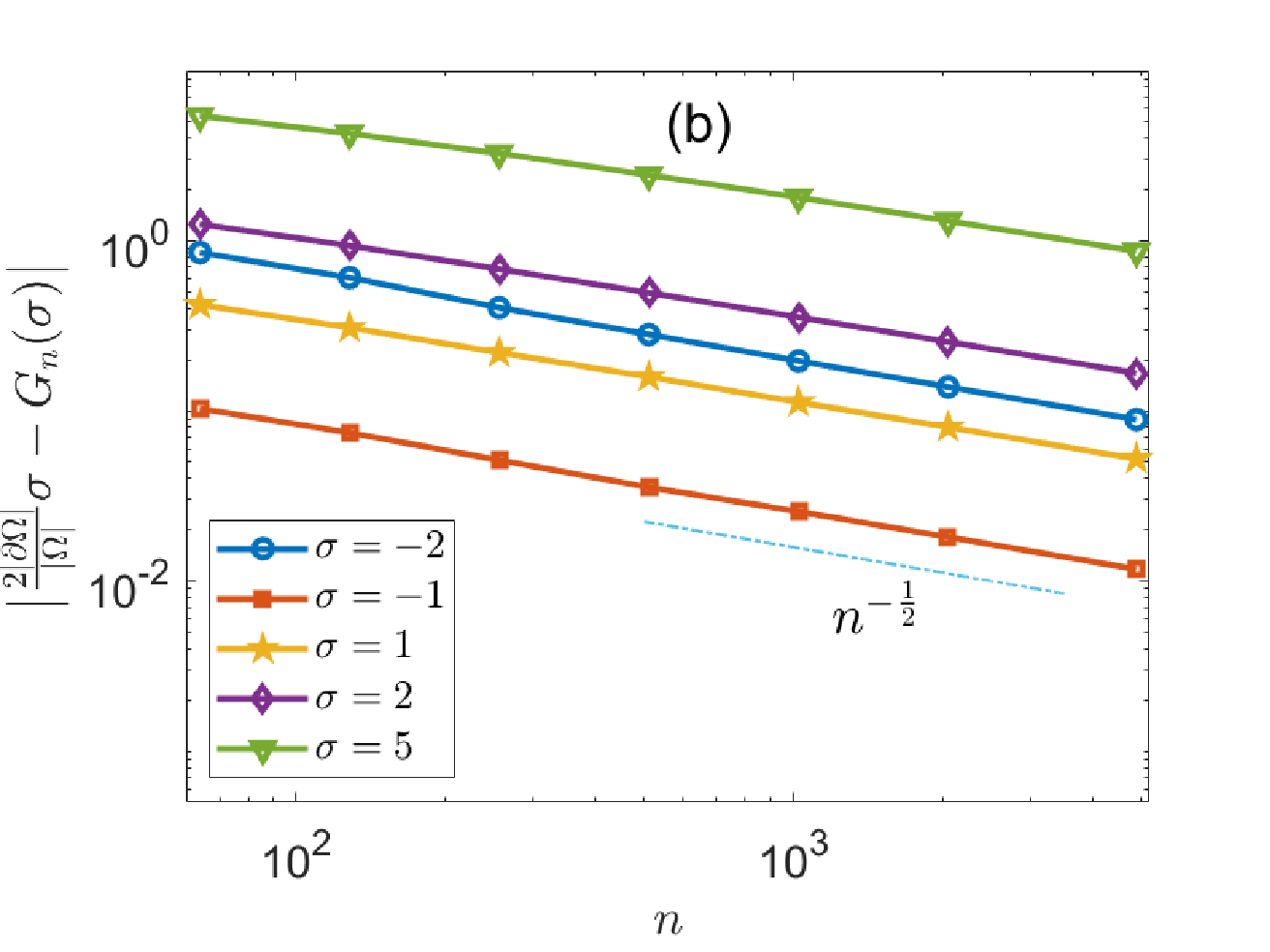}}
  \caption{ Convergence rate of the cumulative averages $G_n(\sigma)$ of \eqref{SO}-\eqref{SORobin} in 2D with
  $V(x,y)\equiv 0$:  (a) with $\sigma=1$ for different domains $\Omega$, and (b) with $\Omega$ taken as a
  quarter-stadium domain for different $\sigma$.}\label{fig:convergence_rate}
 \end{figure}

 Figure \ref{fig:RNgapV1} shows the RtN gaps $g_n(\sigma)$ and $G_n(\sigma)$ of \eqref{SO}-\eqref{SORobin} in 2D with
 $V(x,y)=x^2+y^2$ and $\sigma=1$ for different domains $\Omega$.
 Figure \ref{fig:convergence_rate} shows convergence rates of the
 cumulative averages $G_n(\sigma)$ of \eqref{SO}-\eqref{SORobin} in 2D
 with  $V(x,y)\equiv 0$ and $\sigma=1$ for different domains $\Omega$.

 \begin{figure}[htbp]
 \centering
  \subfigure{ \includegraphics[scale=.32]{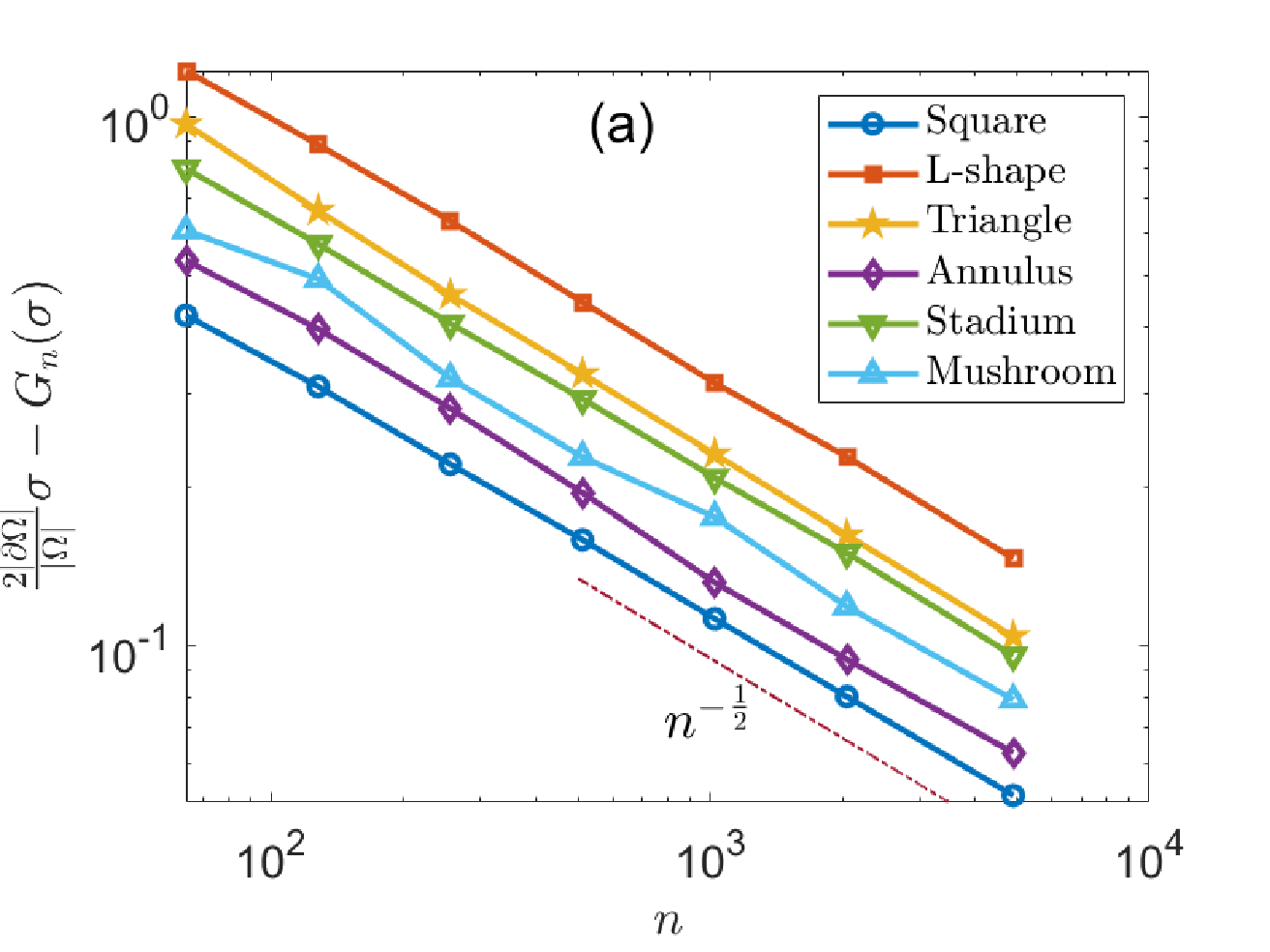}}\hfill
  \subfigure{ \includegraphics[scale=.32]{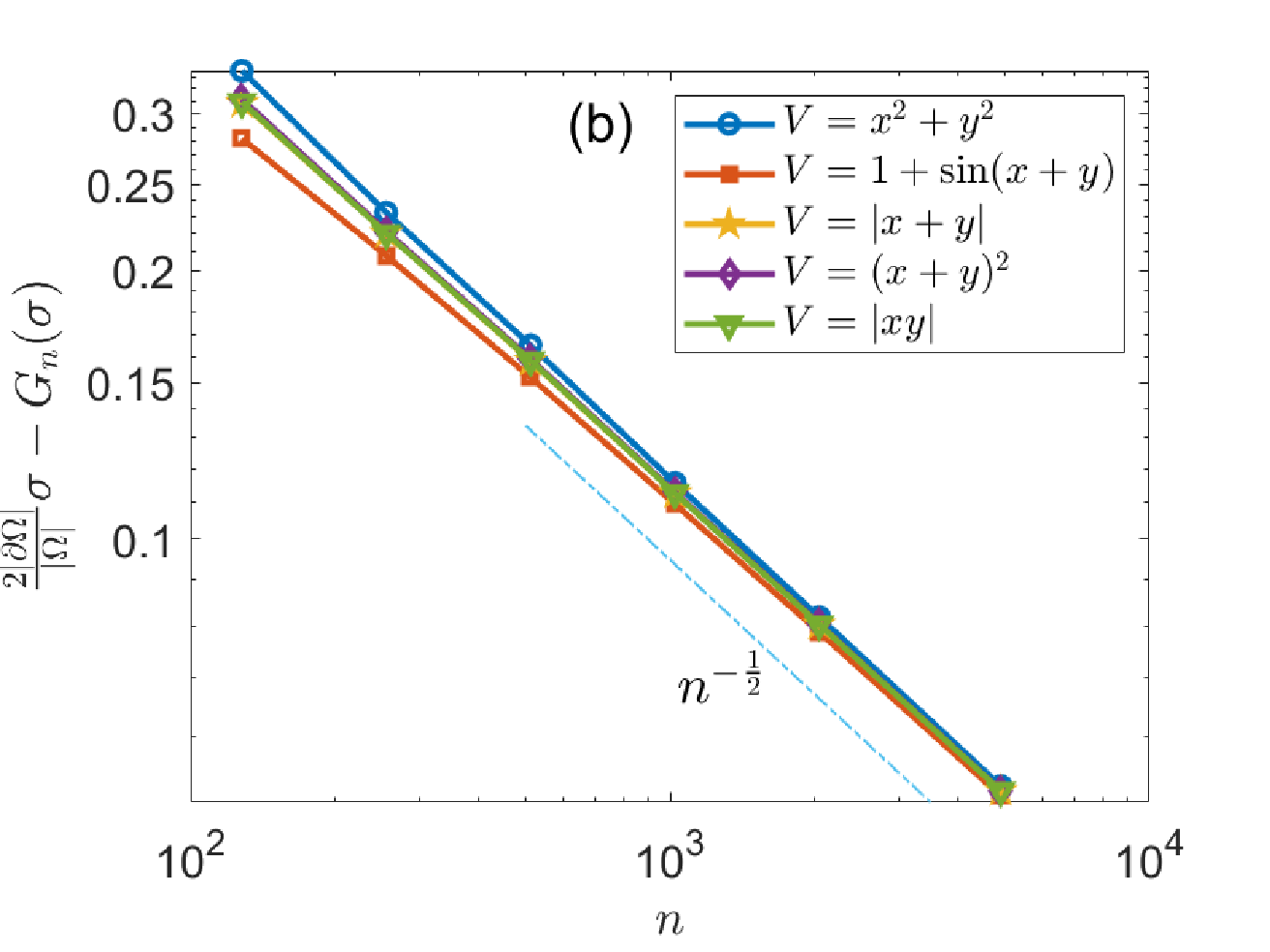}}
  \caption{ Convergence rate of the cumulative averages $G_n(\sigma)$ of \eqref{SO}-\eqref{SORobin} in 2D with
  $\sigma=1$:  (a) with  $V(x,y)=x^2+y^2$ for different domains $\Omega$, and (b) with $\Omega$ taken as a
  quarter-stadium domain for different potentials. }\label{fig:convergence_rateV1}
 \end{figure}

 From Figs. \ref{fig:RNgapV1}-\ref{fig:convergence_rateV1}, and additional results  not shown here for brevity,
we can immediately conclude \eqref{RtNgapdim} when $d=2$.

\subsection{The RtN gaps in 3D}
We take $d=3$ and $V({\bf x})\equiv 0$ in \eqref{SO}-\eqref{SORobin} with two different domains $\Omega$:
one is a cube, i.e. $\Omega=(0,\pi)^3$, and another is a unit ball. The eigenvalues
can be computed semi-analytically from separated transcendental equations.

\begin{figure}[!htbp]
 \centering
  \subfigure{ \includegraphics[scale=.32]{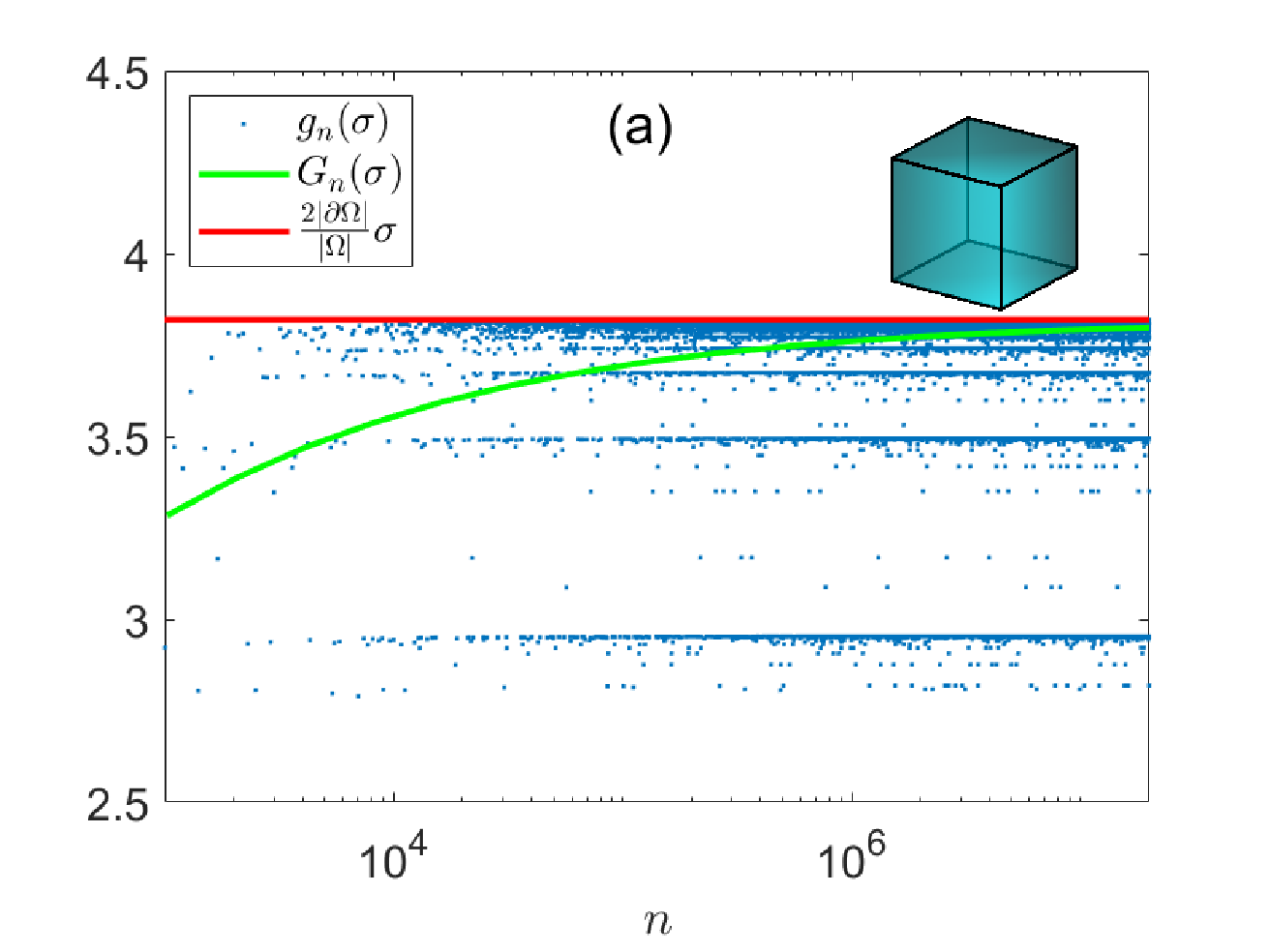}}
  \subfigure{ \includegraphics[scale=.32]{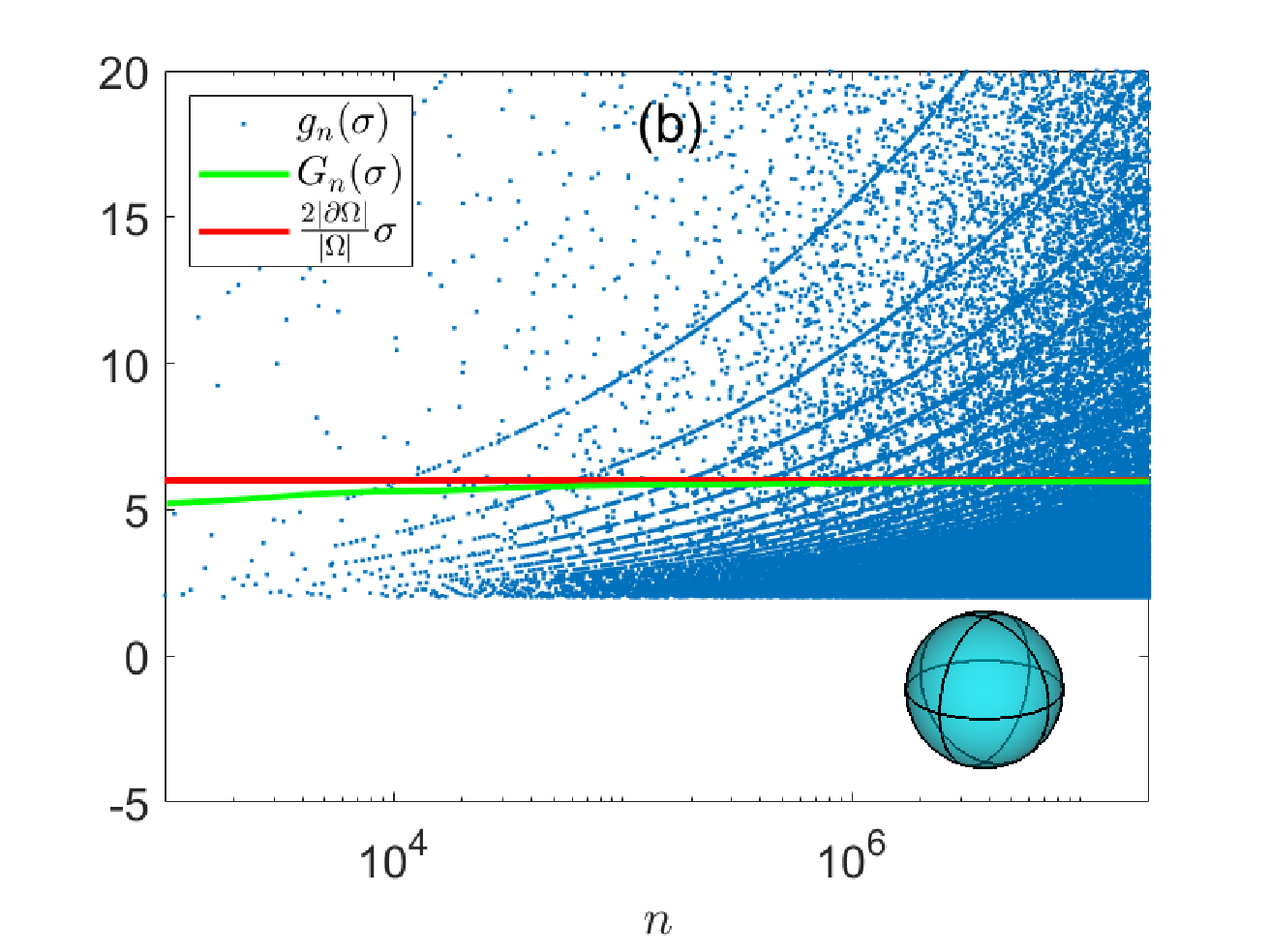}}
  \caption{The RtN gaps of \eqref{SO}-\eqref{SORobin} in 3D with
  $V({\bf x})\equiv 0$ and $\sigma=1$ for different domains $\Omega$:
   (a) the cube $(0,\pi)^3$, and (b) the unit ball.}\label{fig:RNgap3D}
 \end{figure}

\begin{figure}[!htbp]
 \centering
  \subfigure{ \includegraphics[scale=.32]{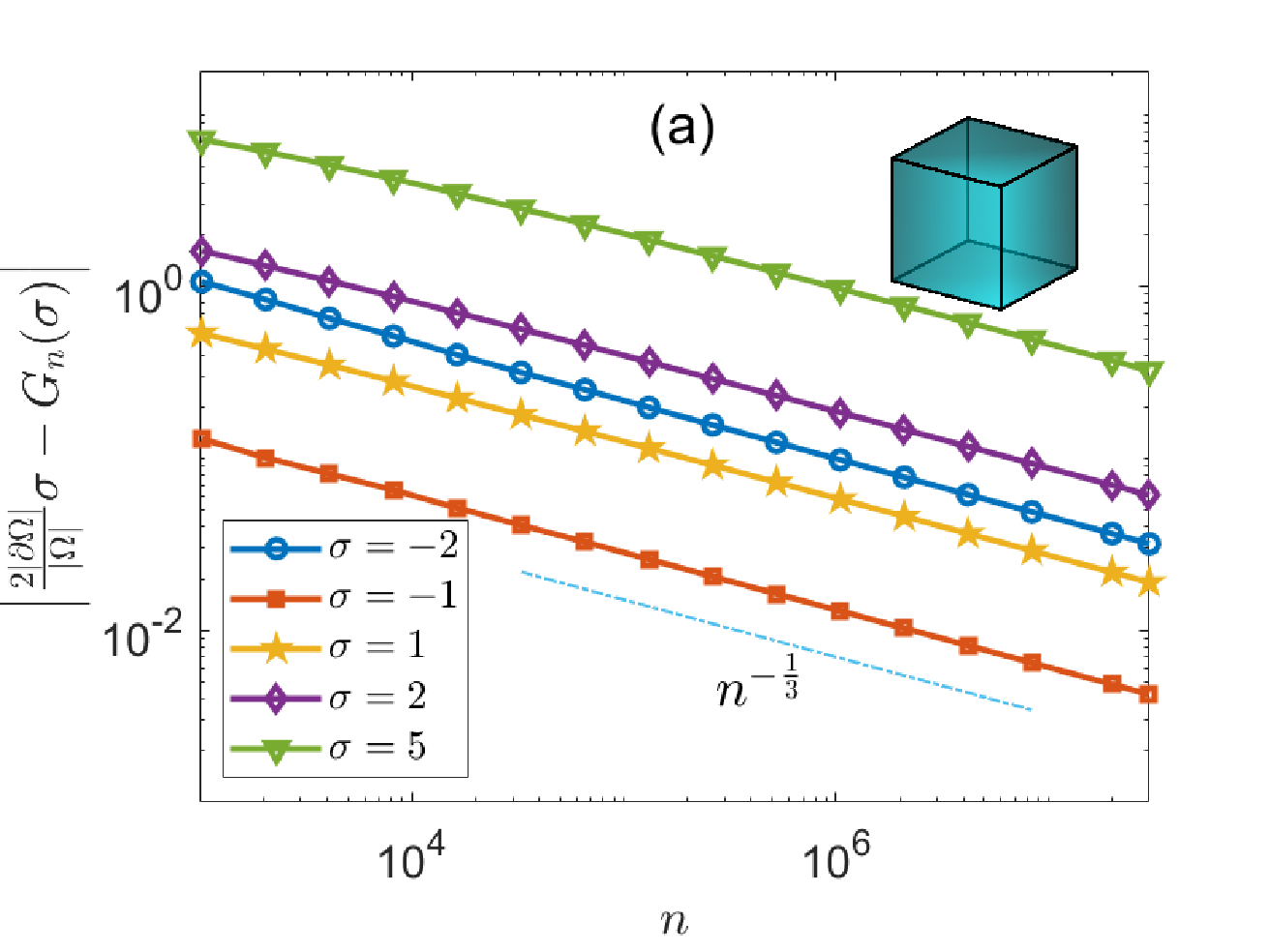}}
  \subfigure{ \includegraphics[scale=.32]{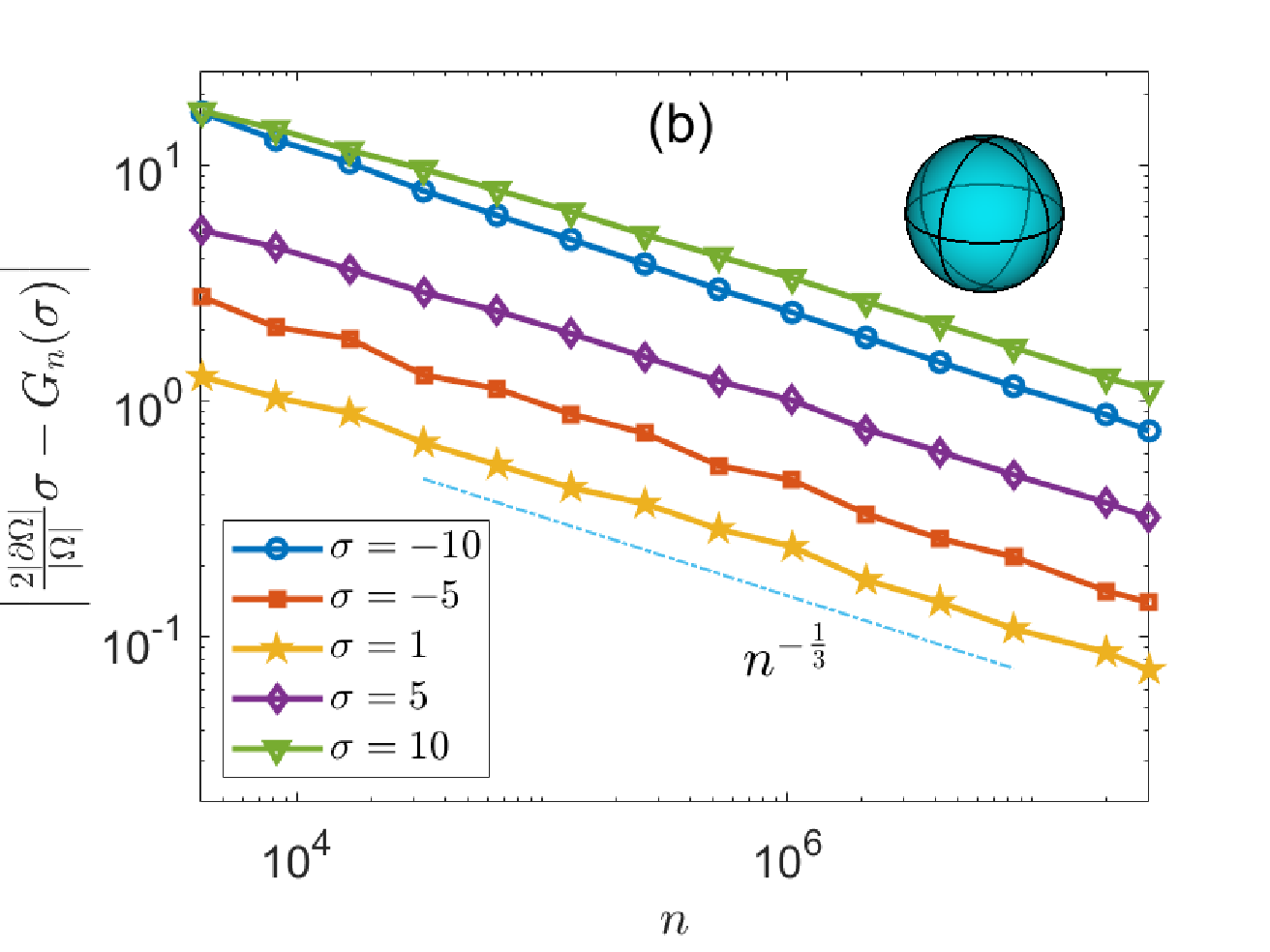}}
  \caption{ Convergence rate of the cumulative averages $G_n(\sigma)$ of \eqref{SO}-\eqref{SORobin} in 3D with
  $V({\bf x})\equiv 0$ for different $\sigma$ and different domains $\Omega$:
  (a) the cube $(0,\pi)^3$, and (b) the unit ball. }\label{fig:conver3D}
 \end{figure}

Figure \ref{fig:RNgap3D} shows the RtN gaps $g_n(\sigma)$ and $G_n(\sigma)$
of \eqref{SO}-\eqref{SORobin} in 3D with $\sigma=1$ for two different domains.
Figure \ref{fig:conver3D} displays convergence rates of the
 cumulative averages $G_n(\sigma)$ of \eqref{SO}-\eqref{SORobin} in 3D
 for different $\sigma$ and  domains $\Omega$.

From Figs. \ref{fig:RNgap3D}-\ref{fig:conver3D}, and additional results  not shown here for brevity,
we can immediately conclude \eqref{RtNgapdim} when $d=3$.

\section{Conclusion}

We proposed and implemented spectral and spectral-element methods for accurately computing many
eigenvalues of the Schr\"odinger operator with Robin boundary condition on simple and complex geometries,
respectively. The proposed framework combines the geometric
flexibility of finite element partitions with the high accuracy of
spectral methods and $p$-version finite element methods,
allowing many reliable eigenvalues to be obtained
under a fixed number of degrees of freedom.
Using thousands of reliable eigenvalues obtained by the proposed numerical methods,
we investigated systematically  the Robin-to-Neumann (RtN) gaps
of Schr\"{o}dinger operator with the Robin boundary condition in one dimension,
two dimensions and three dimensions with different potentials and different geometries.
Based on our extensive numerical results, we formulated
a unified conjecture on the cumulative averages of the RtN gaps of
the Schr\"odinger operator.

\section*{Acknowledgments}
The authors thank Professor Zeev Rudnick for his valuable comments and suggestions, and Professor Zhiguo Yang for his guidance throughout the development of the spectral-element code.
This work was partially supported by the Ministry of Education of Singapore under its AcRF Tier 1 funding grant A-8003584-00-00 (W. Bao) and by the National Natural Science Foundation of China Grant No. 12501554 (F. Huang).

\newpage
\bibliographystyle{plain}
\bibliography{Eigenvalue}

@misc{Shao2026,
  author       = {Li, Y. and Shao, S.},
  title        = {Asymptotic recovery in {Fourier} spectral methods for the {Schr{\"o}dinger} equation
                  with point singularities},
  howpublished = {arXiv:2606.01718},
  year         = {2026},
  doi          = {10.48550/arXiv.2606.01718}
}

@article{lin2025highly,
  author       = {Lin, S. and Li, H. and Yang, Z.},
  title        = {Highly efficient {Gauss}'s law-preserving spectral algorithms for {Maxwell}'s
                  double-curl source and eigenvalue problems based on eigen-decomposition},
  journal      = {CSIAM Trans. Appl. Math.},
  volume       = {6},
  number       = {3},
  pages        = {489--526},
  year         = {2025},
  doi          = {10.4208/csiam-am.SO-2024-0042}
}

@article{han2013class,
  author       = {Han, J. and Yang, Y.},
  title        = {A class of spectral element methods and its a priori/a posteriori error estimates for
                  2nd-order elliptic eigenvalue problems},
  journal      = {Abstr. Appl. Anal.},
  volume       = {2013},
  number       = {1},
  pages        = {262010},
  year         = {2013},
  doi          = {10.1155/2013/262010}
}

@article{yang2016seamless,
  author       = {Yang, Z. and Wang, L.-L. and Rong, Z. and Wang, B. and Zhang, B.},
  title        = {Seamless integration of global {Dirichlet-to-Neumann} boundary condition and spectral
                  elements for transformation electromagnetics},
  journal      = {Comput. Methods Appl. Mech. Engrg.},
  volume       = {301},
  pages        = {137--163},
  year         = {2016},
  doi          = {10.1016/j.cma.2015.12.020}
}

@article{sauter2010hp,
  author       = {Sauter, S. A.},
  title        = {{$hp$}-finite elements for elliptic eigenvalue problems: error estimates which are
                  explicit with respect to $\lambda$, $h$, and $p$},
  journal      = {SIAM J. Numer. Anal.},
  volume       = {48},
  number       = {1},
  pages        = {95--108},
  year         = {2010},
  doi          = {10.1137/070702515}
}

@article{patera1984spectral,
  author       = {Patera, A. T.},
  title        = {A spectral element method for fluid dynamics: laminar flow in a channel expansion},
  journal      = {J. Comput. Phys.},
  volume       = {54},
  number       = {3},
  pages        = {468--488},
  year         = {1984},
  doi          = {10.1016/0021-9991(84)90128-1}
}

@article{rudnick2021rectangles,
  author       = {Rudnick, Z. and Wigman, I.},
  title        = {The {Robin} problem on rectangles},
  journal      = {J. Math. Phys.},
  volume       = {62},
  number       = {11},
  pages        = {113503},
  year         = {2021},
  doi          = {10.1063/5.0061763}
}

@article{rudnick2022hemisphere,
  author       = {Rudnick, Z. and Wigman, I.},
  title        = {On the {Robin} spectrum for the hemisphere},
  journal      = {Ann. Math. Qu{\'e}.},
  volume       = {46},
  number       = {1},
  pages        = {121--137},
  year         = {2022},
  doi          = {10.1007/s40316-021-00155-9}
}

@article{rudnick2022triangle,
  author       = {Rudnick, Z. and Wigman, I.},
  title        = {On the {Robin} spectrum for the equilateral triangle},
  journal      = {J. Phys. A: Math. Theor.},
  volume       = {55},
  number       = {25},
  pages        = {254004},
  year         = {2022},
  doi          = {10.1088/1751-8121/ac6f9b}
}

@article{vanmaele1995operator,
  author       = {Vanmaele, M. and Van Keer, R.},
  title        = {An operator method for a numerical quadrature finite element approximation for a class
                  of second-order elliptic eigenvalue problems in composite structures},
  journal      = {RAIRO Mod{\'e}l. Math. Anal. Num{\'e}r.},
  volume       = {29},
  number       = {3},
  pages        = {339--365},
  year         = {1995},
  doi          = {10.1051/m2an/1995290303391},
  url          = {https://www.numdam.org/item/M2AN_1995__29_3_339_0/}
}

@article{hernandez2003neumann,
  author       = {Hern{\'a}ndez, E. and Rodr{\'i}guez, R.},
  title        = {Finite element approximation of spectral problems with {Neumann} boundary conditions
                  on curved domains},
  journal      = {Math. Comp.},
  volume       = {72},
  number       = {243},
  pages        = {1099--1115},
  year         = {2003},
  doi          = {10.1090/S0025-5718-02-01467-9}
}

@article{lopezyela2017finite,
  author       = {L{\'o}pez-Yela, A. and P{\'e}rez-Pardo, J. M.},
  title        = {Finite element method to solve the spectral problem for arbitrary self-adjoint
                  extensions of the {Laplace--Beltrami} operator on manifolds with a boundary},
  journal      = {J. Comput. Phys.},
  volume       = {347},
  pages        = {235--260},
  year         = {2017},
  doi          = {10.1016/j.jcp.2017.06.043}
}

@article{bogli2022eigenvalues,
  author       = {B{\"o}gli, S. and Kennedy, J. B. and Lang, R.},
  title        = {On the eigenvalues of the {Robin} {Laplacian} with a complex parameter},
  journal      = {Anal. Math. Phys.},
  volume       = {12},
  number       = {1},
  pages        = {39},
  year         = {2022},
  doi          = {10.1007/s13324-022-00646-0}
}

@incollection{babuska1991finite,
  author       = {Babu{\v{s}}ka, I. and Osborn, J.},
  title        = {Eigenvalue problems},
  editor       = {Ciarlet, P. G. and Lions, J.-L.},
  booktitle    = {Finite Element Methods (Part 1)},
  series       = {Handbook of Numerical Analysis},
  volume       = {2},
  pages        = {641--787},
  publisher    = {North-Holland},
  address      = {Amsterdam},
  year         = {1991},
  doi          = {10.1016/S1570-8659(05)80042-0}
}

@article{ognibene2025asymptotics,
  author       = {Ognibene, R.},
  title        = {On asymptotics of {Robin} eigenvalues in the {Dirichlet} limit},
  journal      = {Comm. Partial Differential Equations},
  volume       = {50},
  number       = {9},
  pages        = {1174--1210},
  year         = {2025},
  doi          = {10.1080/03605302.2025.2536098}
}

@book{teschl2014mathematical2,
  author       = {Teschl, G.},
  title        = {Mathematical Methods in Quantum Mechanics: With Applications to {Schr{\"o}dinger}
                  Operators},
  edition      = {Second},
  series       = {Graduate Studies in Mathematics},
  volume       = {157},
  publisher    = {American Mathematical Society},
  address      = {Providence, RI},
  year         = {2014},
  doi          = {10.1090/gsm/157}
}

@book{brenner2008mathematical,
  author       = {Brenner, S. C. and Scott, L. R.},
  title        = {The Mathematical Theory of Finite Element Methods},
  edition      = {Third},
  series       = {Texts in Applied Mathematics},
  volume       = {15},
  publisher    = {Springer},
  address      = {New York},
  year         = {2008},
  doi          = {10.1007/978-0-387-75934-0}
}

@article{larson2000posteriori,
  author       = {Larson, M. G.},
  title        = {A posteriori and a priori error analysis for finite element approximations of
                  self-adjoint elliptic eigenvalue problems},
  journal      = {SIAM J. Numer. Anal.},
  volume       = {38},
  number       = {2},
  pages        = {608--625},
  year         = {2000},
  doi          = {10.1137/S0036142997320164}
}

@article{dai2008convergence,
  author       = {Dai, X. and Xu, J. and Zhou, A.},
  title        = {Convergence and optimal complexity of adaptive finite element eigenvalue computations},
  journal      = {Numer. Math.},
  volume       = {110},
  number       = {3},
  pages        = {313--355},
  year         = {2008},
  doi          = {10.1007/s00211-008-0169-3}
}

@article{carstensen2012adaptive,
  author       = {Carstensen, C. and Gedicke, J.},
  title        = {An adaptive finite element eigenvalue solver of asymptotic quasi-optimal computational
                  complexity},
  journal      = {SIAM J. Numer. Anal.},
  volume       = {50},
  number       = {3},
  pages        = {1029--1057},
  year         = {2012},
  doi          = {10.1137/090769430}
}

@article{carstensen2014guaranteed,
  author       = {Carstensen, C. and Gedicke, J.},
  title        = {Guaranteed lower bounds for eigenvalues},
  journal      = {Math. Comp.},
  volume       = {83},
  number       = {290},
  pages        = {2605--2629},
  year         = {2014},
  doi          = {10.1090/S0025-5718-2014-02833-0}
}

@article{chen2011finite,
  author       = {Chen, H. and He, L. and Zhou, A.},
  title        = {Finite element approximations of nonlinear eigenvalue problems in quantum physics},
  journal      = {Comput. Methods Appl. Mech. Engrg.},
  volume       = {200},
  number       = {21--22},
  pages        = {1846--1865},
  year         = {2011},
  doi          = {10.1016/j.cma.2011.02.008}
}

@article{yang2021eigenfunction,
  author       = {Yang, B. and Zhou, A.},
  title        = {Eigenfunction behavior and adaptive finite element approximations of nonlinear
                  eigenvalue problems in quantum physics},
  journal      = {ESAIM Math. Model. Numer. Anal.},
  volume       = {55},
  number       = {1},
  pages        = {209--227},
  year         = {2021},
  doi          = {10.1051/m2an/2020078}
}

@article{hashemi2022least,
  author       = {Hashemi, B. and Nakatsukasa, Y.},
  title        = {Least-squares spectral methods for {ODE} eigenvalue problems},
  journal      = {SIAM J. Sci. Comput.},
  volume       = {44},
  number       = {5},
  pages        = {A3244--A3264},
  year         = {2022},
  doi          = {10.1137/21M1445934}
}

@book{berezin1991schrodinger,
  author       = {Berezin, F. A. and Shubin, M. A.},
  title        = {The {Schr{\"o}dinger} Equation},
  series       = {Mathematics and Its Applications (Soviet Series)},
  volume       = {66},
  publisher    = {Kluwer Academic Publishers},
  address      = {Dordrecht},
  year         = {1991},
  doi          = {10.1007/978-94-011-3154-4}
}

@book{davies1995spectral,
  author       = {Davies, E. B.},
  title        = {Spectral Theory and Differential Operators},
  series       = {Cambridge Studies in Advanced Mathematics},
  volume       = {42},
  publisher    = {Cambridge University Press},
  address      = {Cambridge},
  year         = {1995},
  doi          = {10.1017/CBO9780511623721}
}

@article{bao2019fundamental,
  author       = {Bao, W. and Ruan, X. and Shen, J. and Sheng, C.},
  title        = {Fundamental gaps of the fractional {Schr{\"o}dinger} operator},
  journal      = {Commun. Math. Sci.},
  volume       = {17},
  number       = {2},
  pages        = {447--471},
  year         = {2019},
  doi          = {10.4310/CMS.2019.v17.n2.a7},
  url          = {https://archive.intlpress.com/site/pub/files/_fulltext/journals/cms/2019/0017/0002/CMS-2019-0017-0002-a007.pdf}
}

@book{griffiths2018introduction,
  author       = {Griffiths, D. J. and Schroeter, D. F.},
  title        = {Introduction to Quantum Mechanics},
  edition      = {Third},
  publisher    = {Cambridge University Press},
  address      = {Cambridge},
  year         = {2018},
  doi          = {10.1017/9781316995433}
}

@article{shan2017triangular,
  author       = {Shan, W. and Li, H.},
  title        = {The triangular spectral element method for {Stokes} eigenvalues},
  journal      = {Math. Comp.},
  volume       = {86},
  number       = {308},
  pages        = {2579--2611},
  year         = {2017},
  doi          = {10.1090/mcom/3173}
}

@article{gittins2020courant,
  author       = {Gittins, K. and Helffer, B.},
  title        = {{Courant}-sharp {Robin} eigenvalues for the square: the case with small {Robin}
                  parameter},
  journal      = {Ann. Math. Qu{\'e}.},
  volume       = {44},
  number       = {1},
  pages        = {91--123},
  year         = {2020},
  doi          = {10.1007/s40316-019-00120-7}
}

@book{shen2011spectral,
  author       = {Shen, J. and Tang, T. and Wang, L.-L.},
  title        = {Spectral Methods: Algorithms, Analysis and Applications},
  series       = {Springer Series in Computational Mathematics},
  volume       = {41},
  publisher    = {Springer},
  address      = {Berlin, Heidelberg},
  year         = {2011},
  doi          = {10.1007/978-3-540-71041-7}
}

@book{deville2002high,
  author       = {Deville, M. O. and Fischer, P. F. and Mund, E. H.},
  title        = {High-Order Methods for Incompressible Fluid Flow},
  series       = {Cambridge Monographs on Applied and Computational Mathematics},
  volume       = {9},
  publisher    = {Cambridge University Press},
  address      = {Cambridge},
  year         = {2002},
  doi          = {10.1017/CBO9780511546792}
}

@article{boffi2010finite,
  author       = {Boffi, D.},
  title        = {Finite element approximation of eigenvalue problems},
  journal      = {Acta Numer.},
  volume       = {19},
  pages        = {1--120},
  year         = {2010},
  doi          = {10.1017/S0962492910000012}
}

@article{li2017efficient,
  author       = {Li, H. and Zhang, Z.},
  title        = {Efficient spectral and spectral element methods for eigenvalue problems of
                  {Schr{\"o}dinger} equations with an inverse square potential},
  journal      = {SIAM J. Sci. Comput.},
  volume       = {39},
  number       = {1},
  pages        = {A114--A140},
  year         = {2017},
  doi          = {10.1137/16M1069596}
}

@article{wang2022spectral,
  author       = {Wang, W. and Zhang, Z.},
  title        = {Spectral element methods for eigenvalue problems based on domain decomposition},
  journal      = {SIAM J. Sci. Comput.},
  volume       = {44},
  number       = {2},
  pages        = {A689--A719},
  year         = {2022},
  doi          = {10.1137/20M1345980}
}

@article{rudnick2021differences,
  author       = {Rudnick, Z. and Wigman, I. and Yesha, N.},
  title        = {Differences between {Robin} and {Neumann} eigenvalues},
  journal      = {Comm. Math. Phys.},
  volume       = {388},
  number       = {3},
  pages        = {1603--1635},
  year         = {2021},
  doi          = {10.1007/s00220-021-04248-y}
}

@article{bao2020jacobi,
  author       = {Bao, W. and Chen, L. and Jiang, X. and Ma, Y.},
  title        = {A {Jacobi} spectral method for computing eigenvalue gaps and their distribution
                  statistics of the fractional {Schr{\"o}dinger} operator},
  journal      = {J. Comput. Phys.},
  volume       = {421},
  pages        = {109733},
  year         = {2020},
  doi          = {10.1016/j.jcp.2020.109733}
}

@article{weideman1988eigenvalues,
  author       = {Weideman, J. A. C. and Trefethen, L. N.},
  title        = {The eigenvalues of second-order spectral differentiation matrices},
  journal      = {SIAM J. Numer. Anal.},
  volume       = {25},
  number       = {6},
  pages        = {1279--1298},
  year         = {1988},
  doi          = {10.1137/0725072}
}

@article{zhang2015many,
  author       = {Zhang, Z.},
  title        = {How many numerical eigenvalues can we trust?},
  journal      = {J. Sci. Comput.},
  volume       = {65},
  number       = {2},
  pages        = {455--466},
  year         = {2015},
  doi          = {10.1007/s10915-014-9971-5}
}

@article{yang2011two,
  author       = {Yang, Y. and Bi, H.},
  title        = {Two-grid finite element discretization schemes based on shifted-inverse power method
                  for elliptic eigenvalue problems},
  journal      = {SIAM J. Numer. Anal.},
  volume       = {49},
  number       = {4},
  pages        = {1602--1624},
  year         = {2011},
  doi          = {10.1137/100810241}
}

@article{lin2015multi,
  author       = {Lin, Q. and Xie, H.},
  title        = {A multi-level correction scheme for eigenvalue problems},
  journal      = {Math. Comp.},
  volume       = {84},
  number       = {291},
  pages        = {71--88},
  year         = {2015},
  doi          = {10.1090/S0025-5718-2014-02825-1}
}

@article{hu2011acceleration,
  author       = {Hu, X. and Cheng, X.},
  title        = {Acceleration of a two-grid method for eigenvalue problems},
  journal      = {Math. Comp.},
  volume       = {80},
  number       = {275},
  pages        = {1287--1301},
  year         = {2011},
  doi          = {10.1090/S0025-5718-2011-02458-0}
}

@article{babuvska1989,
  author       = {Babu{\v{s}}ka, I. and Osborn, J. E.},
  title        = {Finite element-{Galerkin} approximation of the eigenvalues and eigenvectors of
                  selfadjoint problems},
  journal      = {Math. Comp.},
  volume       = {52},
  number       = {186},
  pages        = {275--297},
  year         = {1989},
  doi          = {10.1090/S0025-5718-1989-0962210-8}
}

@article{xu2001,
  author       = {Xu, J. and Zhou, A.},
  title        = {A two-grid discretization scheme for eigenvalue problems},
  journal      = {Math. Comp.},
  volume       = {70},
  number       = {233},
  pages        = {17--25},
  year         = {2001},
  doi          = {10.1090/S0025-5718-99-01180-1}
}

@article{guo2024deep,
  author       = {Guo, Y. and Ming, P.},
  title        = {A deep learning method for computing eigenvalues of the fractional {Schr{\"o}dinger}
                  operator},
  journal      = {J. Syst. Sci. Complex.},
  volume       = {37},
  number       = {2},
  pages        = {391--412},
  year         = {2024},
  doi          = {10.1007/s11424-024-3250-9}
}

@article{guo2026generalization,
  author       = {Guo, Y. and Ming, P. and Yu, H.},
  title        = {Generalization error estimates of a machine learning method for solving
                  high-dimensional {Schr{\"o}dinger} eigenvalue problems},
  journal      = {Sci. China Math.},
  volume       = {69},
  number       = {4},
  pages        = {1063--1114},
  year         = {2026},
  doi          = {10.1007/s11425-024-2434-9}
}

@article{polya1961vibrating,
  author       = {P{\'o}lya, G.},
  title        = {On the eigenvalues of vibrating membranes},
  journal      = {Proc. London Math. Soc. (3)},
  volume       = {11},
  number       = {1},
  pages        = {419--433},
  year         = {1961},
  doi          = {10.1112/plms/s3-11.1.419}
}

@article{li1983schrodinger,
  author       = {Li, P. and Yau, S.-T.},
  title        = {On the {Schr{\"o}dinger} equation and the eigenvalue problem},
  journal      = {Comm. Math. Phys.},
  volume       = {88},
  number       = {3},
  pages        = {309--318},
  year         = {1983},
  doi          = {10.1007/BF01213210}
}

@misc{wang2026toward,
  author       = {Wang, Y. and Xie, H.},
  title        = {Toward {P{\'o}lya}'s conjecture: improving the individual {Li--Yau} bound via energy
                  orthogonality},
  howpublished = {arXiv:2609.17219},
  year         = {2026},
  doi          = {10.48550/arXiv.2609.17219}
}

@misc{wang2026improved,
  author       = {Wang, Y. and Xie, H.},
  title        = {Improved {Berezin--Li--Yau} bounds for {Dirichlet} eigenvalues via {Fourier} coupling
                  and frequency concentration},
  howpublished = {arXiv:2611.01741},
  year         = {2026},
  doi          = {10.48550/arXiv.2611.01741}
}

@article{jiang2026polya,
  author       = {Jiang, R. and Lin, F.},
  title        = {{P{\'o}lya}'s conjecture up to {$\epsilon$}-loss and quantitative estimates for the
                  remainder of {Weyl}'s law},
  journal      = {Comm. Pure Appl. Math.},
  eid          = {e70058},
  year         = {2026},
  note         = {Article e70058 (Early View)},
  doi          = {10.1002/cpa.70058}
}

@article{zhai2019lower,
  author       = {Zhai, Q. and Zhang, R.},
  title        = {Lower and upper bounds of {Laplacian} eigenvalue problem by weak {Galerkin} method on
                  triangular meshes},
  journal      = {Discrete Contin. Dyn. Syst. Ser. B},
  volume       = {24},
  number       = {1},
  pages        = {403--413},
  year         = {2019},
  doi          = {10.3934/dcdsb.2018091}
}

@article{carstensen2020skeletal,
  author       = {Carstensen, C. and Zhai, Q. and Zhang, R.},
  title        = {A skeletal finite element method can compute lower eigenvalue bounds},
  journal      = {SIAM J. Numer. Anal.},
  volume       = {58},
  number       = {1},
  pages        = {109--124},
  year         = {2020},
  doi          = {10.1137/18M1212276}
}

@article{filonov2023balls,
  author       = {Filonov, N. and Levitin, M. and Polterovich, I. and Sher, D. A.},
  title        = {{P{\'o}lya}'s conjecture for {Euclidean} balls},
  journal      = {Invent. Math.},
  volume       = {234},
  number       = {1},
  pages        = {129--169},
  year         = {2023},
  doi          = {10.1007/s00222-023-01198-1}
}

@article{filonov2026annuli,
  author       = {Filonov, N. and Levitin, M. and Polterovich, I. and Sher, D. A.},
  title        = {{P{\'o}lya}'s conjecture for {Dirichlet} eigenvalues of annuli},
  journal      = {J. London Math. Soc. (2)},
  volume       = {113},
  number       = {2},
  pages        = {e70425},
  year         = {2026},
  doi          = {10.1112/jlms.70425}
}
\end{document}